\documentclass[12pt,a4paper]{amsart}
\usepackage[utf8]{inputenc}
\usepackage[english]{babel}
\usepackage[T1]{fontenc}
\usepackage{amsmath,amsthm,amsfonts,amssymb, hyperref}
\usepackage{color}
\usepackage[left=2cm,right=2cm,top=2cm,bottom=2cm]{geometry}
\usepackage{cleveref}

\usepackage[pdftex]{graphicx}
\usepackage{tikz-cd}
\usepackage{float}
\usepackage[autostyle]{csquotes}
   \MakeAutoQuote{‘}{’}

\numberwithin{equation}{subsection}
\numberwithin{figure}{subsection}

\newtheorem{theorem}[equation]{Theorem}
\newtheorem{definition}[equation]{Definition}
\newtheorem{proposition}[equation]{Proposition}
\newtheorem{corollary}[equation]{Corollary}
\newtheorem{lemma}[equation]{Lemma}
\newtheorem{remark}[equation]{Remark}

\newtheorem{conjecture}[equation]{Conjecture}
\newtheorem{question}[equation]{Question}

\newtheorem{theorem2}{Theorem}[section]

\newtheorem{proposition2}[theorem2]{Proposition}

\newtheorem{lemma2}[theorem2]{Lemma}
\newtheorem{remark2}[theorem2]{Remark}

\newcommand{\QED}{\hfill $\square$}

\newcommand{\dd}{\mathrm{d}}

\newcommand{\Ham}{\mathrm{Ham}}

\newcommand{\N}{\mathbb{N}}
\newcommand{\Z}{\mathbb{Z}}

\newcommand{\R}{\mathbb{R}}

\newcommand{\U}{\mathcal{U}}
\newcommand{\F}{\mathcal{F}}
\newcommand{\V}{\mathcal{V}}
\newcommand{\GG}{\mathcal{G}}

\allowdisplaybreaks  %allows page breaking
\begin{document}
\title[pb invariant on surfaces]{The Poisson Bracket Invariant on Surfaces}
\author[J. Payette]{Jordan Payette}
\address{Department of Mathematics and Statistics\\
McGill University\\
Montréal, Québec, Canada}
\email{jordan.payette@mail.mcgill.ca}

\begin{abstract}
We study the Poisson bracket invariant, which measures the level of Poisson noncommutativity of a smooth partition of unity, on closed symplectic surfaces. Motivated by a general conjecture of Polterovich \cite{P3} and building on preliminary work of Buhovsky--Tanny \cite{BT}, we prove that for any smooth partition of unity subordinate to an open cover by discs of area at most $c$, and under some localization condition on the cover when the surface is a sphere, then the product of the Poisson bracket invariant with $c$ is bounded from below by a universal constant. Similar results were obtained recently by  Buhovsky--Logunov--Tanny \cite{BLT} for open covers consisting of displaceable sets on all closed surfaces, and their approach was extended by Shi—Lu \cite{SL} to open covers by nondisplaceable discs. We investigate the sharpness of all these results.
\end{abstract}

\maketitle

\tableofcontents

\section{Introduction and results} \label{sec:Introduction}

In the course of his investigation of function theory on symplectic manifolds, Polterovich \cite{P2} introduced the so-called \textit{Poisson bracket invariant} $pb(\F)$ as a quantitative measure of Poisson non-commutativity of functions forming a partition of unity $\F$ on a symplectic manifold. He also explained the relation of this symplectic invariant to operational quantum mechanics, more precisely describing how the invariant appears as a lower bound on the statistical noise of measurements of the collections of quantum observables associated to $\F$ via Berezin-Toeplitz quantization. This relation highlights the importance of establishing lower bounds on the Poisson bracket invariant. \vspace{6pt}

The present paper considers this problem when the symplectic manifolds are closed surfaces equipped with an area form. Our main result, refined by \Cref{cor_pb_thm}, is: \vspace{6pt}

\noindent \textsc{Poisson bracket theorem:}\label{main_thm} There exist constants $C, C' > 0$ such that the following holds. Let $(M, \omega)$ be a closed symplectic surface and $\F$ be a partition of unity subordinate to a finite open cover $\U = \{U_i\}_{i \in I}$ by topological discs of area no larger than $c$ . Let $X = \{x_1, \dots, x_m\} \subset M$ be points at which $\U$ localizes \textit{i.e.} each disc $U_i \in \U$ contains no more than one $x_k \in X$. (Note that $m \ge 1$.)  When $M = S^2$, assume $m \ge 3$. Then
\[ pb(\F) \, c \; \ge \; C \; \mbox{ and } \;  pb(\F) \, \mathrm{Area}(M, \omega)  \, \ge \,  C' \, m \; . \] \vspace{0pt}

\noindent This result is closely related to (and was motivated by) a general question of Polterovich \cite{P3} which in the special case of closed symplectic surfaces amounts to asking if the estimate $pb(\F) \, c \; \ge \; C$ holds for a universal constant $C > 0$ in the ``small-scale regime'' when $\U$ consists of displaceable discs, \textit{i.e.} when $c \le \mathrm{Area}(M, \omega)/2$. A positive answer to Polterovich's question for all closed symplectic surfaces was recently given by Buhovsky--Logunov--Tanny \cite{BLT}. Our main theorem, which holds without restriction on $c$, also confirms Polterovich's conjecture for closed surfaces of genus $g \ge 1$ as well as for covers on the sphere that satisfies our localization conditions. Though our main theorem does not fully recover Buhovsky--Logunov--Tanny's result on the sphere, it is of interest in other respects: \vspace{4pt}

(i)  Our main theorem generalizes Buhovsky--Logunov--Tanny's result to the extent that it shows that the geometric condition of displaceability of the discs in the open cover is not necessary in general and that it can be relaxed to more topological (\textit{i.e.} diffeomorphism-invariant) conditions. \vspace{4pt}

(ii) With one exception, our estimates are proved with constants $C$ and $C'$ that are as large, and usually larger, than those obtained in \cite{BLT}. \vspace{4pt}

(iii) Our main theorem generalizes another result from \cite{BLT} (already in \cite{BT}): Our main theorem holds if $m$ denotes instead the number of (displaceable) ``essential discs'' in $\U$, where $U \in \U$ is essential if $\U \setminus \{U\}$ does not cover $M$, assuming $\U$ admits such a disc (\textit{i.e.} assuming $m \ge 1$). \vspace{4pt}

\noindent The starting point of our method is to notice that this last result holds with the displaceability condition relaxed to a topological condition we name ``confinement''. Topological operations then allow to deduce the main theorem from this special case. Let's turn to a more detailed description of all this.

\subsection{Preliminary notions}\label{sec:preliminaries}

A \emph{symplectic manifold} $(M, \omega)$ is a smooth manifold $M$ endowed with a closed nondegenerate differential 2-form $\omega$, \textit{i.e.} $\dd \omega = 0$ and the bundle map $\omega_{\flat} : TM \to T^{\ast}M : X \mapsto X \lrcorner \, \omega$ is an isomorphism. The \emph{symplectic gradient} of a smooth function $H$ on $M$ is the vector field $X_H$ on $M$ defined as $-\dd H = X_H \lrcorner \, \omega$. The \emph{Poisson bracket} associated to $(M, \omega)$ is the bilinear map
\begin{align}
\notag \{ - , - \} \, :& \; C^{\infty}(M) \times C^{\infty}(M) \to C^{\infty}(M) \\
\notag & (G,H) \mapsto \{G,H\} := -\omega(X_G, X_H) \, . 
\end{align}
\noindent A \emph{symplectomorphism} of $(M, \omega)$ is a diffeomorphism $\phi$ of $M$ that preserves $\omega$ in the sense that $\phi^*\omega = \omega$. Hence symplectomorphisms preserve the volume form $\wedge^n\omega/n!$ A diffeomorphism $\phi : M \to M$ is \emph{Hamiltonian} if there exists a smooth time-dependent function $h_t$ on $M$ (compactly supported in the interior of $M$), defined for $t \in (- \epsilon, 1+\epsilon)$ for some $\epsilon > 0$, such that $\phi = \phi^1_h$ where $\phi^t_h$ solves the Cauchy problem $\phi^0_h = Id$ and $\frac{d}{dt} \phi^t_h(x) = X_h \circ \phi^t_h(x)$. The set of Hamiltonian diffeomorphisms, denoted $\Ham(M, \omega)$, is a subgroup of the identity component $\mathrm{Symp}_0(M, \omega)$ of the group of symplectomorphisms (see for instance \cite{Ba,P1,PR}).\vspace{12pt}

We say that a subset $X \subset (M, \omega)$ is \emph{displaceable} if there exists $\phi \in \Ham(M, \omega)$ which \emph{displaces} $X$ (from its closure $\overline{X}$), \textit{i.e.} $\phi(X) \cap \overline{X} = \varnothing$. Another invariant associated to a set $X \subset (M, \omega)$ is its (\emph{Hofer}) \emph{displacement energy}, denoted $e_H(X)$. Namely,  first define the \emph{Hofer norm} or \emph{energy} of a diffeomorphism $\phi \in \Ham(M, \omega)$ as
\begin{align}
\notag \| \phi \|_H := \underset{h : \phi^1_h = \phi}{\inf} \;   \int_0^1  \left( \max_{x \in M} h_t(x) - \min_{x \in M} h_t(x)\right)dt \, , 
\end{align}
then set 
\[ e_H(X) := \inf \, \{ \, \| \phi \|_H \, : \, \phi \in \Ham(M, \omega), \, \phi(X) \cap \overline{X} = \varnothing \,  \} \, ,\]
with the convention that $e_H(X) = + \infty$ if $X$ is nondisplaceable. \vspace{12pt}

\begin{remark}\label{rmk:def_displaceability}
Slightly different notions of a `displaceable subset' appear in the literature. The original notion due to Hofer \cite{H} is more restrictive than ours: A set $X \subset (M, \omega)$ is `displaceable in Hofer's sense' if its closure is displaceable in our sense. A more permissive definition than ours, used \emph{e.g.} in \cite{P1}, considers as `displaceable' any $X$ such that $\phi(X) \cap X = \varnothing$ for some $\phi \in \Ham(M, \omega)$. Perhaps the most suitable notion of a `displaceable subset' is one having finite `displacement energy', with the latter concept defined as in \cite{U}. Our definition, borrowed from \cite{PR}, strikes a balance between generality and simplicity of use.
\end{remark}\vspace{12pt}

\subsection{Poisson bracket invariant} \label{sec:pbdefinition}

Given a smooth function $f : M \to \R$, its \emph{support} is $\mathrm{supp}(f) := \overline{\{ x \in M \, | \, f(x) \neq 0\}}$. We say that $f$ is \emph{supported in an open set $U$} if $\mathrm{supp}(f) \subseteq U$. A (\emph{smooth}) \emph{positive collection} is a locally finite collection of functions $\F = \{ f_i \}_{i \in I}$ such that each $f_i \ge 0$ and 
\[ S_{\F}(x) := \sum_{i \in I} f_i(x) \ge 1 \quad \mbox{ for all } x \in M \, .\] A (\emph{smooth}) \emph{partition of unity} is a positive collection $\F$ such that $S_{\F} = 1$. A positive collection $\F$ is said to be \emph{subordinate to a locally finite open cover} $\U = \{U_i\}_{i \in I}$, denoted $\F \prec \U$, if for each $i \in I$, $f_i$ is supported in $U_i$. We denote by $\| - \| : C^{0}(M;\R) \to [0, +\infty)$ the supremum norm on real-valued continuous functions on $M$. \vspace{12pt}

Following Polterovich \cite{P2}, define the \emph{Poisson bracket invariant of a positive collection $\F = \{f_i\}_{i\in I}$} as
\[ pb(\F) :=  \underset{a,b \in [-1,1]^I}{\mathrm{sup}} \, \left\| \left\{ \sum_{i \in I} a_i f_i \, , \, \sum_{j \in I} b_j f_j \right\} \right\| \; , \]
(this supremum is achieved and finite when $M$ is compact) and define the \emph{Poisson bracket invariant of an open cover $\U$} as the infimum of $pb(\F)$ over positive collections $\F$ subordinate to $\U$,
\[ pb(\U) := \inf_{\F \prec \, \U} \, pb(\F) \; .  \]
These quantities are symplectic invariants to the extent that given a symplectomorphism $\phi$ of $(M, \omega)$, denoting $\phi_* \F = \{(\phi^{-1})^*f_i\}_{i \in I}$ and $\phi_*\U = \{\phi(U_i)\}_{i \in I}$, we have $pb(\phi_*\F) = pb(\F)$ and $pb(\phi_*\U) = pb(\U)$.\vspace{12pt}

In this paper, we shall mainly consider a related invariant introduced by Buhovsky and Tanny \cite{BT}. Define the \emph{Poisson bracket function of a positive collection $\F = \{f_i\}_{i\in I}$} as
\[ P_{\F} : M \to [0, \infty) : x \mapsto \sum_{i,j \in I} \left| \{f_i, f_j \}(x) \right| \, . \]
It has been established that there exists a constant $0 < c(n) \le 1$ (depending only on the dimension $\mathrm{dim} \, M = 2n$) such that for any positive collection,
\[  c(n) \, \|P_{\F}\| \le pb(\F) \le \|P_{\F}\| \, . \]
The upper bound is a straightforward application of the triangle inequality. The lower bound was first established in dimension $2$ in \cite{BT}, then in every dimension in \cite{BLT} and was greatly improved in \cite{GlT}. Therefore, the problem of finding lower bounds on $pb(\F)$ is equivalent to that of finding lower bounds on $\|P_{\F}\|$. Of course, lower bounds on the $L^{\infty}$-norm of $P_{\F}$ follow from lower bounds on its $L^1$-norm of $P_{\F}$; As in \cite{BLT}, we are going to establish lower bounds on $\|P_{\F}\|_{L^1}$. We shall see that this approach is essentially optimal in dimension two, in the sense that lower bounds on $\|P_{\F}\|_{L^1}$ conversely follow from lower bounds on $\|P_{\F}\|$ (see \Cref{Linfty_Lone}). \vspace{12pt}

\subsection{Poisson bracket conjectures}\label{sec:pbconjecture}

Building upon previous work (\emph{e.g.} \cite{EP, EPZ, P2}; see \cref{sec:history}), Polterovich \cite{P3} asked whether the following statement is true:

\begin{conjecture}[Polterovich's Poisson bracket conjecture]\label{conj} There exists a constant $C > 0$ depending only on $(M, \omega)$ such that for any open cover $\U$ consisting of displaceable sets, setting $e_H(\U) = \mathrm{sup}_{U_i \in \U} \, e_H(U_i)$,
\begin{equation} \label{PoC}  pb(\U)  \ge C \, e_H(\U)^{-1} \; . \end{equation}
\end{conjecture} \vspace{0pt}

\noindent A related and \textit{a priori} simpler conjecture is:\vspace{0pt}

\begin{conjecture}[Weak Poisson bracket conjecture]\label{weakconj} There exists a constant $C > 0$ depending only on $(M, \omega)$ such that for any open cover $\U$ consisting of displaceable sets,
\begin{equation} \label{PoC2}  pb(\U)  \ge C \, \mathrm{Vol}(M, \omega)^{-1/n} \; . \end{equation}
\end{conjecture}\vspace{0pt}

These conjectures formulate precise versions of the following plausible general reciprocity relation: \emph{The more localized an open cover is, the bigger its $pb$ invariant should be}. In Polterovich's conjecture, the localization of an open cover is measured by the displacement energy of its constituents. The significance of the second conjecture is to claim that $pb$ is bounded away from $0$ -- \emph{by a constant independent of the cover} -- on covers consisting of displaceable sets. We note that the conjectures have been formulated so that the constant $C$ would depend on $\omega$ only up to symplectomorphisms and multiplication by  positive scalars, hence in fact only up to weak deformation equivalence, see \cite[\S 13.2]{MS}.  \vspace{12pt}

Both conjectures have been more or less explicitly discussed in \cite{BT, BLT, EP, EPZ, I, LPa, Pa, P2, P3, P4, PR, Se, SL}. Albeit we shall give a more detailed description of these works in \cref{sec:history}, we give here a brief motivation for these conjectures. Note that $pb(\U)$ may vanish if $\U$ does not consist of displaceable sets: Think of a partition of unity subordinate to a cover consisting of only two open sets. On the other hand, Polterovich \cite{P3} has established that whenever the open cover is formed by displaceable sets, then the above inequalities hold with $C$ a positive $\U$-dependent quantity.\vspace{12pt}

The conjectures are more tractable in dimension two, since in that case both displaceability and the Poisson bracket can be understood in more elementary terms. In that dimension, both conjectures have recently been proved in \cite{BLT} after a preliminary breakthrough in \cite{BT}. We shall now turn to explaining how the confirmation of the conjectures for special covers in \cite{BT} can be combined with topological arguments to obtain lower bounds on the pb invariant of open covers satisfying other, more topological notions of localization that both generalize and refine that given by the displacement energy.

\subsection{Main concepts and results} \label{sec:Results}

In order to state our results precisely, we need a few more definitions.\vspace{0pt}

\begin{definition}
Let $(M, \omega)$ be a closed symplectic surface and $\U = \{U_1, \dots, U_N \}$ be an open cover. The \textbf{capacity of $\U$} is
\[ c(\U) \, := \,  \underset{1 \le i \le N}{\mathrm{max}} \, \mathrm{Area}(U_i, \omega) \; . \]\vspace{-12pt}
\end{definition}

\begin{definition}
Let $\U = \{U_i\}_{i \in I}$ be a locally finite open cover of a smooth manifold $M$. A connected subset $X \subset M$ is \textbf{confined} (\textbf{with respect to $\U$}) if some connected component of $\partial \overline{X}$ is not contained in any single $U \in \U$.
\end{definition}

\begin{definition}
Let $\U = \{U_i\}_{i \in I}$ be a locally finite open cover of a smooth manifold $M$. Given $x \in M$, denote
\[ \U_x \, := \,  \{  \, U_i \in \U \, : \, x \in U_i \, \}  \]
and define the \textbf{star of $x$} (\textbf{in $\U$}) to be the region
\[ St(x) = St(x;\U) \, := \,  \bigcup_{U_i \in \U_x} \, U_i \, \subseteq \, M \; . \]
We denote $\Gamma_{\U}$ the set of points $x \in M$ for which $St(x)$ is confined.
\end{definition}

\begin{definition}
Let $\U = \{U_i\}_{i \in I}$ be a locally finite open cover of a symplectic surface $M$. Given a positive collection $\F \prec \U$, the \textbf{star function} $C_{\F} : M \to [0, \infty)$ is the measurable function defined as
\[ C_{\F}(x) :=  \int_{St(x)} \sum_{i \in I} \sum_{j : U_j \in \U_x} |\{f_i, f_j\}| \, \omega \quad (x \in M) \, . \]
\end{definition}

\begin{remark}\label{rmk:def_star}
Morally, $St(x)$ is confined if it does not spread in $M$ so much as for its boundary components to be included in single sets of $\U$: It is ``confined to stay near $x$''. We also note that, since $\U$ is locally finite, each $C_{\F}$ attains its maximum over any compact set in $M$.
\end{remark}\vspace{0pt}

\begin{definition}
Let $\U = \{U_i\}_{i \in I}$ be a locally finite open cover of a smooth manifold $M$. We say that $\U$ is \textbf{localized at points $x_1, \dots, x_m$} if each $U \in \U$ contains at most one of these points. For $m \in \N$, we say that $\U$ is \textbf{localized in $m$ points} or is \textbf{$m$-localized} if there are $m$ points $x_1, \dots, x_m \in M$ at which $\U$ is localized.
\end{definition} \vspace{0pt}

\begin{remark}
Any open cover is $1$-localized. An open cover $\U$ is localized at the points $x_1, \dots, x_m$ if and only if for all $1 \le i,j \le m$, $x_i \in St(x_j)$ implies $x_i = x_j$. For $m > 1$, a $m$-localized cover is in particular $(m-1)$-localized, merely by forgetting about one $x_k$.
\end{remark}\vspace{12pt}

The next two notions are borrowed from \cite{BT,BLT}:\vspace{0pt}

\begin{definition}
Let $\U = \{U_i\}_{i \in I}$ be a locally finite open cover of a smooth surface $M$. A disc $U \in \U$ is \textbf{essential} (\textbf{to $\U$}) if $\U \setminus \{U\}$ is not a cover of $M$. Equivalently, there exists $x \in U$ such that $\U_x = \{U\}$, \textit{i.e.} such that $St(x) = U$.
\end{definition} \vspace{0pt}

\begin{definition}
Let $M$ be a smooth surface. A locally finite open cover $\U = \{U_i\}_{i \in I}$ on $M$ is said to be \textbf{in general position} if the sets $U_i$ have smooth boundaries, if the boundaries intersect transversally \emph{i.e.} $\partial U_i \pitchfork \partial U_j$ for all $i \neq j$, and if $\partial U_i \cap \partial U_j \cap \partial U_k = \varnothing$ for every triple $(i,j,k)$ of distinct indices.
\end{definition} \vspace{0pt}

\begin{remark}\label{rmk:gen_position}
A suitably general setting for our forthcoming results would consider locally finite open covers $\U = \{U_i\}_{i \in I}$ of smooth surfaces by relatively compact topological discs. By standard smooth approximation and transversality arguments, we can approximate each $U_i \in \U$ by an open set $U'_i \subset U_i$ in such a way that $\U' := \{U'_i\}_{i \in I}$ is in general position. Of course, $e_H(\U') \le e_H(\U)$ and $c(\U') \le c(\U)$. The perturbation of $\U$ into $\U'$ can be made in such a way that discs remain discs, while preserving the essentiality of discs, the confinement of stars and the localization at given points (which we may arrange to lie in the complement of $\cup_{i \in I} \partial \overline{U'}_i$). Moreover, given $\F \prec \U$, since each $f_i$ has compact support in the open set $U_i$, we may arrange to have $\F \prec \U'$. Hence, there is no real loss of generality in assuming from the outset that our open covers are in such general positions, and we shall henceforth do so.
\end{remark}\vspace{12pt}

We now come to our results. The starting point of this paper is the following inequality established in \cite{BT,BLT}, although stated here in a slightly more general form: We allow for the surface to be open, for the open cover to be locally finite, for the functions to form a positive collection rather than a partition of unity, and we replaced a displaceability condition by a confinement condition. \vspace{12pt} 

\begin{theorem}[``Confined essential disc'' inequality]\label{conf_ess}
Let $(M, \omega)$ be a (possibly open) symplectic surface, $\U = \{U_i\}_{i \in I}$ be a locally finite open cover consisting of discs in general position and $\F = \{f_i\}_{i \in I} \prec \U$ be a positive collection. Assume that the set $\mathcal{J}_c(\U) \subset \U$ of confined essential discs is nonempty. Then
\begin{align}
\notag \int_{U_j} \sum_{i \in I} |\{f_i, f_j\}| \, \omega \, &\ge \, 1 \; \mbox{ for all } \, U_j \in \mathcal{J}_c(\U)\,.
\end{align}
\end{theorem}\vspace{12pt}

\begin{corollary}[\cite{BT, BLT}]\label{cor_conf_ess}
Under the same assumption as in \Cref{conf_ess}, we have
\[ \underset{U_i \in \U}{\mathrm{sup}} \, \int_{U_i} P_{\F} \, \omega \ge \underset{U_i \in \mathcal{J}_c(\U)}{\mathrm{inf}} \, \int_{U_i} P_{\F} \, \omega \ge 1 \; \mbox{ and } \; \int_M P_{\F} \, \omega \ge |\mathcal{J}_c(\U)| \; . \]
\end{corollary} \vspace{12pt}

\textit{Proof}. The first string of inequalities follows from $P_{\F} \ge \sum_{i \in I} |\{f_i, f_j\}|$ and \Cref{conf_ess}. The second inequality follows similarly, observing that $|\{f_i, f_j\}| = 0$ outside $U_j$, so that
\begin{align}
\notag  \int_M P_{\F} \, \omega &\ge \int_M \sum_{i \in I}\sum_{j : U_j \in \mathcal{J}_c(\U)} |\{f_i, f_j\}| \, \omega \\
\notag &= \sum_{j : U_j \in \mathcal{J}_c(\U)} \int_{U_j} \sum_{i \in I} |\{f_i, f_j\}| \, \omega \\
\notag &\ge \sum_{j : U_j \in \mathcal{J}_c(\U)} 1 = |\mathcal{J}_c(\U)| \, . 
\end{align}

\QED
\vspace{12pt}

A central idea in our approach is to lift open covers (by discs) and positive collections along (ramified) symplectic covering maps, since this procedure tends to simplify the topology of the stars $St(x)$ to the point that they become confined, while preserving (or increasing in a controlled way) quantities such as $\|P_{\F}\|$ and $c(\U)$. This leads to a proof of the following \textit{star inequality}, which in a sense semi-locates the Poisson non-commutativity on the surface and which could also be understood as a generalization of the essential inequality.\vspace{12pt}

\begin{theorem}[Star inequality] \label{star_ineq}
Let $(M, \omega)$ be a (possibly open) symplectic surface, $\U = \{U_i\}_{i \in I}$ be a locally finite open cover consisting of open discs in general position and $\F \prec \U$ be a positive collection. 
\begin{itemize}
 \item If $M \neq S^2$ and $x \in M$, then $C_{\F}(x)  \, \ge  \, 1$.
\item If $M=S^2$ and $x \in \Gamma_{\U}$, then $C_{\F}(x)  \, \ge  \, 1$.
\item If $M = S^2$ and $\U$ is $m$-localized with $m \ge 3$, say at points $X = \{x_1, x_2, x_3, \dots, x_m\}$, then for any $x \in S^2$,
\begin{align}\label{eq_star_ineq2}
C_{\F}(x)  \, &\ge  \, C = \begin{cases} 1/3 & \mbox{ if $m=3$}, \\ 1/2 &  \mbox{ if $m\ge 4$}, \\ \;\, 1 & \mbox{ if $x$ is close to $X$. }
\end{cases}
\end{align}
Here, `$x$ close to $X$' means that $x$ belongs to the same connected component of $M \setminus \cup_{i \in I} \partial \overline{U_i}$ as some $x_k \in X$.
 \end{itemize}
\end{theorem}\vspace{0pt}

\begin{remark}\label{rmk:superfluous}
When $M=S^2$, the assumptions of confinement on $St(x)$ or of $3$-localization on $\U$ are relevant (although they might be relaxable, see \Cref{question:star}). Indeed, consider the open cover $\U = \{U_n, U_s\}$ on $S^2 = \{(x,y,z) \in \R^3 \, : \, x^2+y^2+z^2 = 1\}$, where $U_{n} = \{z  > -1/2\}$ and $U_s = \{ z < 1/2 \}$. It is only $2$-localized and the possible star sets are $U_n$, $U_s$ and $S^2$, none of which being confined. However, for any partition of unity $\F \prec \U$, one has $P_{\F} \equiv 0$. 
\end{remark}\vspace{0pt}

Next, consider a nonzero positive measure $\mu$ on $M$. Let $\mu(h)$ denote the integral of $h : M \to \R$ over $M$ with respect to $\mu$. We have:

\begin{theorem}[Abstract Poisson bracket theorem] \label{pb_thm}
Let $(M, \omega)$ be a closed symplectic surface, $\U = \{U_1, \dots, U_N\}$ be an open cover of $M$ by discs in general position and $\F \prec \U$ be a positive collection. Let $\mu$ be a nonzero finite positive Borel measure on $M$ and set $\mu(\U) = \mathrm{max}_{1 \le i \le N} \, \mu(U_i)$. Then
\begin{align}
\notag \int_{M} P_{\F} \, \omega  \, &\ge \;   \frac{\mu(C_{\F})}{\mu(\U)} \; .
\end{align}
\end{theorem}\vspace{12pt}

Combining \Cref{star_ineq} and \Cref{pb_thm} with $\mu$ being either the measure induced by $\omega$ or a uniform discrete measure on the set of points at which the open cover localizes, and using the forthcoming \Cref{Linfty_Lone}, we deduce our main result:

\begin{theorem}[Poisson bracket theorem] \label{cor_pb_thm}
Let $(M, \omega)$ be a closed surface, $\U = \{U_i\}_{i \in I}$  be a finite open cover by discs localized at the points $X = \{x_1, \dots, x_m\}$, and $\F \prec \U$ be a positive collection. If $M = S^2$ and $m \le 2$, assume that $X \subset \Gamma_{\U}$. Then
\[ \int_{M} P_{\F} \, \omega  \, \ge \;  m \; . \] 
Moreover, denoting $ \mathbb{E}[P_{\F}] := (\mathrm{Area}(M, \omega))^{-1} \, \int_{M} P_{\F} \, \omega  $, the estimate
\[ \mathbb{E}[P_{\F}] \, c(\U) \ge C  \quad \mbox{\textit{i.e.}} \quad \int_{M} P_{\F} \, \omega  \, \ge \;  \frac{C \, \mathrm{Area}(M, \omega)}{c(\U)} \,  , \]
holds with $C = 1$ if $M$ has genus $g \ne 1$ or if $\Gamma_{\U} = M$, while for $M = S^2$, it holds with $C=1/2$ if $m \ge 4$ and with $C=1/3$ if $m=3$. This estimate in turn implies
\[ \|P_{\F} \| \, c(\U) \, \ge \, C \, , \; \mbox{ and in fact } \;   \underset{U_i \in \U}{\mathrm{max}} \, \int_{U_i} P_{\F} \, \omega \, \ge \, C \, . \]
\end{theorem}\vspace{0pt}

\textit{Proof}. Without loss of generality, we assume that $\U$ is in general position. We shall use the lower bounds on $C_{\F}$ given by \Cref{star_ineq}.

Given that $\U$ is localized at the points $x_1, \dots, x_m$, let $\mu = \sum_{k=1}^m \delta_{x_k}$, where $\delta_{x_k}$ denotes the Dirac measure supported at $x_k$. Since we have $C_{\F}(x) \ge 1$ for all $x$ sufficiently close to some point $x_k$, we obtain $\mu(C_{\F}) \ge \mu(M) = m$ and $\mathrm{max}_{1 \le i \le N} \, \mu(U_i) = 1$ since each $U_i \in \U$ contains at most one $x_k$. Hence \Cref{pb_thm} implies
\[ \int_M P_{\F} \, \omega \, \ge \, m  \, . \]

Next, let $\mu$ in \Cref{pb_thm} be the measure determined by the symplectic form $\omega$ on $M$. Then $\mathrm{max}_{1 \le i \le N} \, \mu(U_i) = c(\U)$ and $\mu(C_{\F}) \ge C \, \mathrm{Area}(M, \omega)$ holds  with $C = 1$ if $M$ has genus $g \ge 1$ or if $\Gamma_{\U} = M$, with $C = 1/2$ if $M=S^2$ and $m \ge 4$, and with $C = 1/3$ if $M = S^2$ and $m=3$, since $C_{\F} \ge C$ over $M$ in all those situations. Hence \Cref{pb_thm} implies
\[ \int_M P_{\F} \, \omega \, \ge \,  C \,  \frac{\mathrm{Area}(M, \omega)}{c(\U)} \, . \]
It therefore follows that $\|P_{\F}\| \, c(\U) \, \ge \, C$ in those different cases. Since the $m$-localization of a cover is a diffeomorphism-invariant condition,  view of \Cref{Linfty_Lone}, we also have
\[  \underset{U_i \in \U}{\mathrm{max}} \, \int_{U_i} P_{\F} \, \omega \, \ge \, C \, .   \]

\QED
\vspace{12pt}

\begin{remark}\label{rmk:cor_pb_thm}
For any finite open cover $\U$ by discs on any closed surface $(M, \omega)$ and $\F \prec \U$, we also have the estimate
\[ \|P_{\F}\| \, c(\U) \, \ge \, \mathbb{E}[P_{\F}] \, c(\U) \, \ge \, \frac{\mathrm{Area}(\Gamma_{\U}, \omega)}{\mathrm{Area}(M, \omega)} \, . \]
Despite being of interest, we decided not to include this estimate in the main statement, since the lower bound is not invariant under diffeomorphisms when it is smaller than $1$ (hence we cannot apply \Cref{Linfty_Lone} to deduce a lower bound on $\mathrm{max}_{U_i \in \U} \int_{U_i} P_{\F} \, \omega$).
\end{remark} \vspace{12pt}

\subsection{Poisson bracket conjectures in dimension two} \label{sec:discussion_conj}

The above results imply the correctness of conjectures \ref{conj} and \ref{weakconj} in the case of closed symplectic surfaces of genus $g \ge 1$:

\begin{corollary} \label{cor_pb_conj}
Let $(M, \omega)$ be a closed symplectic surface of genus $g \ge 1$ equipped with an open cover $\U = \{U_1, \dots, U_N\}$ by displaceable sets. Then given any positive collection $\F \prec \U$,
\[  \int_{M} P_{\F} \, \omega \, \ge \, \frac{\mathrm{Area}(M, \omega)}{e_H(\U)} \ge 2 \, .  \]
Consequently, for every $\F \prec \U$,
\[ \|P_{\F}\| \, \mathrm{Area}(M, \omega) \, \ge \, 2 \; \mbox{ and } \; \|P_{\F}\| \, e_H(\U) \, \ge \, 1 \, .   \]
\end{corollary} \vspace{12pt}

\textit{Proof}. Fix $\epsilon > 0$. For each $i = 1, \dots, N$, the closed set $S_i := \mathrm{supp}(f_i)$ lies inside $U_i$ and is thus displaceable, with $e_H(S_i) \le e_H(U_i)$. We now rely on the well known characterization of displaceability in dimension two (see \Cref{sec:disp_two} for a proof): Given a closed symplectic surface $(M, \omega)$, a connected closed subset $X \subset M$ is displaceable if and only if it is contained in an embedded closed disc $X'$ of area less than half that of $M$, in which case it is possible to find $X'$ with area arbitrarily close to the displacement energy of $X$. Consequently, each set $S_i$ is contained in an embedded closed disc $U'_i$ that satisfies $\mathrm{Area}(U'_i, \omega) < e_H(U_i) + \epsilon/2 < (1/2)\mathrm{Area}(M, \omega) + \epsilon$. Since $g \ge 1$, \Cref{cor_pb_conj} then follows from \Cref{cor_pb_thm} applied to $\U'$ and from taking the limit as $\epsilon$ goes to zero.

\QED
\vspace{12pt}

\begin{remark}
Our methods confirm Polterovich's conjecture for those open covers $\U$ on $M=S^2$ by displaceable sets that refine open covers $\U'$ by discs localized in $m \ge 3$ points. Indeed, we note that if $\U''$ is another open cover by discs refining $\U'$ and refined by $\U$, then $\U''$ is also localized in $m$ points. Hence, as in the proof above, we may assume $c(\U') < e_H(\U) + \epsilon$. Applying \Cref{cor_pb_thm} to $\U'$ and letting $\epsilon$ go to zero, we get $\|P_{\F}\| e_H(\U) \ge 1/3$ and $\|P_{\F}\| \mathrm{Area}(M, \omega) \ge 3$.
\end{remark} \vspace{12pt}

\subsection{Comparison with the literature} \label{sec:comparison}

Buhovsky--Logunov--Tanny \cite{BLT} proved that for every closed surface (including the sphere), the following lower bound is valid for positive collections $\F$ with displaceable supports,
\begin{equation}\label{BLT_lowerbound} \int_M P_{\F} \, \omega \ge \frac{\mathrm{Area}(M, \omega)}{2e_H(\U)} \quad \mbox{\textit{i.e.}} \quad \mathbb{E}[P_{\F}] \, e_H(\U) \, \ge \, \frac{1}{2} \, ,  \end{equation}
thereby proving the full Poisson bracket conjecture on closed surfaces. For surfaces of genus $g \ge 1$, by lifting the data along a degree-$2$ covering map, it is easily seen that Buhovsky--Logunov--Tanny's result holds verbatim for every open cover by discs, thereby recovering part of \Cref{cor_pb_thm}, albeit with a lower bound half as big as the one we achieved.\vspace{6pt}

For $M = S^2$ and $\U$ an open cover by displaceable discs which is localized in exactly $3$ points, we only proved $\mathbb{E}[P_{\F}] \, e_H(\U) \ge 1/3$, but we achieve the same estimate $\mathbb{E}[P_{\F}] \, e_H(\U) \ge 1/2$ if the covers $4$-localized or satisfies $\mathrm{Area}(\Gamma_{\U}, \omega) \ge \mathrm{Area}(M, \omega)/2$. In fact, for open covers for which $\Gamma_{\U} = M$, we obtain the stronger estimate $\mathbb{E}[P_{\F}] \, e_H(\U) \ge 1$.\vspace{6pt}

Interestingly, Shi and Lu \cite{SL} recently extended the methods of \cite{BLT} to prove that given an open cover $\U$ consisting of finitely many (possibly nondisplaceable) discs such that no two discs cover $M$ (which is a void condition in genus $g \ge 1$), then any positive collection $\F \prec \U$ satisfies $\int_M P_{\F} \, \omega \, \ge \, 2 $. For covers having a large disc \textit{i.e.} with $c(\U) > \mathrm{Area}(M, \omega)/2$, this lower bound implies  $\mathbb{E}[P_{\F}] \, e_H(\U) \ge 1$, which is twice as much as that assured by \Cref{cor_pb_thm}. 
\vspace{0pt}

\subsection{Sharpness of the estimates} \label{subsec:sharpness}

The previous discussion raises the question of the sharpness of all these estimates. Our first observation is that for topologically defined families of open covers, estimating the supremum norm of $P_{\F}$ via its $L^1$-norm is an essentially optimal approach:

\begin{proposition}\label{Linfty_Lone}
Let $(M, \omega)$ be a closed surface. Let $\mathfrak{U}$ be a family of finite open covers on $M$ by discs in general position that is invariant under the action of $\mathrm{Diff}_0(M)$ (the identity component of diffeomorphism group of $M$) given by $(\phi \, , \,  \U = \{U_i\}_{i \in I}) \mapsto  \phi^*\U = \{\phi^{-1}U_i\}_{i \in I}$. Given constants $C, C' > 0$, the following two statements are equivalent:
\begin{enumerate}
\item For all $\U  \in \mathfrak{U}$ and all positive collections $\F \prec \U$,
\[ \| P_{\F} \| \, \mathrm{Area}(M, \omega) \, \ge \, C \; \mbox{ and } \; \|P_{\F} \| \, c(\U) \, \ge \, C' \, ;   \]

\item For all $\U  \in \mathfrak{U}$ and all positive collections $\F \prec \U$,
\[ \int_M P_{\F} \, \omega \, \ge \, C \; \mbox{ and } \; \underset{U_i \in \U}{\mathrm{max}} \, \int_{U_i} P_{\F} \, \omega \, \ge \, C' \, .   \]
\end{enumerate}
\end{proposition}\vspace{0pt}

\begin{remark}\label{rmk:Linfty_Lone} 
The family of finite open covers by discs in general position is invariant under the above action of $\mathrm{Diff}_0(M)$, as are the collection of covers that have a confined star set and the collections (indexed by $m \in \N$) of $m$-localized covers.
\end{remark}\vspace{0pt}

Our next observation is that the ``confined essential disc'' estimate, and hence our best star estimates, are sharp:

\begin{proposition}\label{prop:sharp_conf_ess_disc}
For any symplectic surface $(M, \omega)$ and $\epsilon > 0$, there exist a locally finite open cover by discs $\U = \{U_i\}_{i \in I}$ with $U_1 \in \U$ confined and essential and a positive collection $\F = \{f_i\}_{i \in I} \prec \U$ such that
\[ \int_{U_1} \sum_{i \in I} |\{f_i, f_1\}| \, \omega \, \le \, 1 + \epsilon \, .\]
\end{proposition}

Despite this result, the lower bounds in \Cref{cor_pb_thm} can be improved for partitions of unity in the presence of confined essential discs:

\begin{proposition}\label{prop:lowerbound2}
In the context of \Cref{conf_ess}, suppose that $\F$ is a partition of unity. Then
\[ \underset{U_i \in \U}{\mathrm{max}} \, \int_{U_i} P_{\F} \, \omega \ge \underset{U_i \in \mathcal{J}_c(\U)}{\mathrm{min}} \, \int_{U_i} P_{\F} \, \omega \ge 2 \; . \]
Moreover, if there are $J \ge 1$ disjoint confined essential discs, then
\[ \int_M P_{\F} \, \omega \ge 2J \, . \]
\end{proposition} \vspace{12pt}

In view of this result and of Shi--Lu's estimate, it might be that all lower bounds in \Cref{cor_pb_thm} could be multiplied by a factor $2$. Since \Cref{conf_ess} is a sharp result, improving the lower bounds in \Cref{pb_thm} seems a difficult task using our methods. It would be interesting to see if our techniques could be mixed with those of \cite{BLT, SL} to give a better and more complete understanding of the $pb$ invariant on surfaces. With this in mind, we ask:

\begin{question}\label{question:star}
Given a finite open cover $\U = \{U_i\}_{i=1}^N$ by discs on $S^2$ such that no two discs in $\U$ suffice to cover $S^2$ and a positive collection $\F \prec \U$, is  $C_{\F} \ge 1$ everywhere on $S^2$? (Note that $C_{\F}$ is independent of the choice of symplectic form.)
\end{question} \vspace{0pt}

\subsection{Structure of the paper} \label{sec:structure}

\Cref{sec:proofs} gathers the proofs of our main theorems. We first deduce \Cref{pb_thm} from \Cref{star_ineq} in subsection \ref{sec:proof_pb_thm} using a simple integration argument. Buhovsky--Tanny's ``confined disc inequality'', \Cref{conf_ess}, is proved in the next subsection \ref{sec:proof_conf_ess}. We proceed in subsection \ref{sec:star_genus} to proving \Cref{star_ineq} for all surfaces except the sphere, by lifting the data $(\U, \F)$ to the universal cover along a symplectic map. Subsection \ref{sec:star_sphere} is devoted to the proof of \Cref{star_ineq} in the case of the sphere, first under the assumption of a confined star, then under the assumption of localization in $3$ or more points. The essence of the proof of this last case is the same as that of other surfaces, but we now need to lift the data on the sphere to data on the torus along appropriate ramified symplectic covering maps (of degree $3$ when $m=3$ and of degree $2$ when $m \ge 4$).

In \cref{sec:sharpness}, we prove and discuss our sharpness results from subsection \ref{subsec:sharpness}, in order of appearance. \vspace{6pt}

The paper concludes with two appendices.

\Cref{sec:disp_two} presents a proof of the characterization of displaceable (closed) sets in two dimensions. Although this appears to be a rather well-known result to the experts, we were unable to find a proof of it in the literature and so decided to include a proof of it in this paper.

Finally, \Cref{sec:history} is a short account of our reading of the history of the Poisson bracket invariant and of the Poisson bracket conjectures. It thereby includes some complementary information on the Poisson bracket conjectures in dimension $2$ established in the works \cite{BLT,SL}. We also use this last section to briefly mention how the ideas of the present paper could be used to approach the Poisson bracket conjectures in higher dimensions.

\vspace{12pt}

\textbf{Acknowledgements.} This work has been undertaken as part of the author's PhD thesis at the Université de Montréal. I thank my advisor François Lalonde for encouragement over the years, and also for his invitation at the CIRM in Luminy in June 2015 when he was the Jean-Morlet Chair, on which occasion he introduced me to the Poisson bracket conjectures and to the work of L. Polterovich on the subject. I thank the CIRM for their hospitality during my stay, of which I keep splendid memories. Very special thanks go to Jean Lagacé for his suggestion of a ``mean-to-max'' argument for directly proving a corollary of the Poisson bracket conjectures. Albeit this work has not been included in the final version of this article, his suggestion served as an important guide towards the averaging technique used to prove \Cref{pb_thm}. I warmly thank Dominique Rathel-Fournier for discussions about the characterization of displaceability in two dimensions. The current version of the paper was written and corrected during my postdoctoral positions at Tel Aviv University, under the supervision of Lev Buhovsky and Leonid Polterovich, and at McGill University, under the supervision of Brent Pym; I am deeply grateful to them for their guidance and advice and to these institutions for their stimulating research environments. I also thank the anonymous referee for many useful suggestions that helped improve the presentation and the content of the paper. Over the course of this project, I was supported by the NSERC's Alexander Graham Bell Doctoral Scholarship, by the ISM Excellence Scholarship, by the ERC Starting Grant 757585 and by the FRQNT Postdoctoral Research Scholarship 304098. \vspace{12pt}

%%%%%%%%%%%%%%%%%%%%%%%%%%%%%%%%%%%%555

\section{Proofs of the main results}\label{sec:proofs}

\subsection{Abstract Poisson bracket theorem}\label{sec:proof_pb_thm}

We present here how integration of the star function leads to lower bounds on $\int_M P_{\F} \, \omega$. \vspace{12pt}

\textit{Proof of \Cref{pb_thm}}. We recall that
\[ C_{\F}(x) := \int_{St(x)} \sum_{i=1}^N \sum_{j : U_j \in \U_x} |\{f_i, f_j\}| \, \omega  \quad (x \in M) \, .\]
We now let $\mu$ be a nonzero positive finite (Borel) measure on the closed surface $M$. We shall write $d\mu_x$ to denote the density of $\mu$ at $x \in M$, and similarly write $\omega_y$ to denote the value of the 2-form $\omega$ at $y \in M$.  For a measurable function $h : M \to \R$, let $\mu(h) := \int_M h(x) d\mu_x$. Given a measurable set $S \subset M$, we denote by $\chi_S : M \to \{0,1\}$ its characteristic function. We note that $\mu(U_i) = \mu(\chi_{U_i})$. We also use the shorthand $\mu(\U) = {\mathrm{max}}_{U_i \in \U} \, \mu(U_i)$. Recall that if $y \not \in St(x)$, then for all $1 \le i \le N$ and for all $j$ such that $U_j \in \U_x$, we have $\{f_i, f_j\}(y) = 0$. We compute
\begin{align}
\notag \mu(C_{\F}) &=  \int_{x \in M} \int_{y \in St(x)} \sum_{i=1}^N \sum_{j : U_j \in \U_x} |\{f_i, f_j\}(y)| \, \omega_y \; d\mu_x \\
\notag &= \int_{x \in M} \int_{y \in M} \sum_{i=1}^N \sum_{j : U_j \in \U_x} |\{f_i, f_j\}(y)| \, \omega_y \; d\mu_x \\
\notag &= \sum_{i=1}^N \sum_{j=1}^N \int_{x \in M} \int_{y \in M} \chi_{U_j}(x) |\{f_i, f_j\}(y)| \, \omega_y \; d\mu_x \\
\notag &= \sum_{i=1}^N \sum_{j=1}^N \left( \int_{x \in M} \chi_{U_j}(x) \, d\mu_x \right)  \left( \int_{y \in M} |\{f_i, f_j\}(y)| \, \omega_y \right) \\
\notag &= \sum_{i=1}^N \sum_{j=1}^N \; \mu(U_j)  \; \int_{y \in M} |\{f_i, f_j\}(y)| \, \omega_y \\
\notag & \le \, \mu(\U) \; \int_{y \in M}\sum_{i=1}^N \sum_{j=1}^N \, |\{f_i, f_j\}(y)| \, \omega_y \, .
\end{align}
\noindent This proves the estimate:
\[ \int_M P_{\F} \, \omega \, \ge \, \frac{\mu(C_{\F})}{\mu(\U)} = \underset{1 \le i \le N}{\mathrm{min}} \, \frac{\mu(C_{\F})}{\mu(U_i)}. \]

\QED
\vspace{12pt}

\subsection{``Confined essential disc'' inequality}\label{sec:proof_conf_ess}

We now prove \Cref{conf_ess}. Although the statement is more general and the result stronger than the corresponding results from \cite{BT, BLT}, the proof is essentially the same. \vspace{12pt}

\textit{Proof of \Cref{conf_ess}}. Without loss of generality, suppose $U_1 \in \mathcal{J}_c(\U)$. The idea of the proof consists in using the flow generated by the Hamiltonian vector field associated to $f_1$ in order to dynamically parametrize (a portion of) $U_1$ via ``energy-time'' coordinates, to use Fubini theorem to express the above double integrals as iterated integrals in ``time'' and ``energy'', and finally to use the identity $|\{f_i, f_1\}| = |\mathrm{d}f_i(X_{f_1})|$ to understand the ``time integrals'' as measures of the total oscillation of $f_i$ along each integral curve of $X_{f_1}$.

Let's fix a Riemannian metric on $M$, hence allowing to take the gradient vector field of the smooth function $f_1$. Since $U_1$ is essential, there exists $x \in U_1$ such that no other $U_i \in \U$ contains $x$; Consequently, $f_1(x) = S_{\F}(x) \ge 1$. Also, $f_1(y) = 0$ for $y \in M \setminus U_1$, hence  $(0,1) \subset f_1(M)$ by the intermediate value theorem. By Sard's theorem, the set of regular values $I \subset (0,1)$ of $f_1$ has full measure $1$. Moreover, since $f_1$ is continuous with compact support, $I$ is open and it is therefore a union of disjoint open intervals $I = \cup_{\alpha \in A} \, I_{\alpha}$. For each $\alpha \in A$, choose $s_{\alpha} \in I_{\alpha}$; We wish to pick $p_{\alpha} \in f_1^{-1}(s_{\alpha})$ appropriately, which necessitates a detour.

Let $s$ be any regular value of $f_1$ in $(0,1)$. The regular level-set $f_1^{-1}(s)$ therefore consists of finitely many disjoint smoothly embedded circles $\{C_k\}_{k=1}^K$ included in the open disc $U_1$. By the Jordan-Schoenflies theorem, $U_1 \setminus C_k$ consists of an open disc $D_k$ and of an open annulus $A_k$ such that $\partial \overline{A}_k = C_k \cup \partial U_1$. We claim that as least one $D_k$ contains $x$. Otherwise, $x$ would be contained in the intersection $W := \cap_{k=1}^K A_k$. This intersection is connected: Indeed, $C_l \subset D_k$ implies $D_l \subset D_k$ and inclusions of the discs $D_k$ into one another give a partial order on the discs. So there are finitely many maximal discs, each disjoint from one another, and the complement of their union is the intersection $W$. Contracting each maximal disc to a point within itself, we conclude that $W$ is homotopy equivalent to $U_1$ with a finite number of punctures, which is indeed connected. Consequently, there exists a continuous path $\gamma$ in $W$ from $x$ to $\partial U_1$; By the intermediate value theorem applied to $f_1 \circ \gamma$, we deduce that there exists $y \in W$ such that $f_1(y) = s$, in contradiction to the construction of $W$. Now pick any $D_k \ni x$; We claim that $C_k = \partial D_k$ intersects the complement of any disc $U_i$ in the open cover other than $U_1$. Otherwise, we would have $C_k \subset U_i$ for some $i \neq 1$, hence $C_k$ would bound a closed disc $D'$ in $U_i$ whose interior would necessarily coincide with one of the two connected components of $M \setminus C_k$, one of which being $D_k$. By definition of $x$, we have $x \not \in U_i$, hence $x \not \in  D'$. Thus $D' = M \setminus D_k$ and $U_i \supset D' \supset \overline{A_k} \supset \partial U_1$, in contradiction to the assumption that $U_1$ is confined.

For each $\alpha \in A$, pick $p_{\alpha}$ in a circle of $f_1^{-1}(s_{\alpha})$ winding around $x$. Following the gradient flow line of $f_1$ through $p_{\alpha}$, we get an embedding $\gamma_{\alpha} : I_{\alpha} \to U_1$ such that $\gamma_{\alpha}(s_{\alpha}) = p_{\alpha}$ and $f_1(\gamma_{\alpha}(s)) = s$. We define $\gamma : I \to U_1$ to be the embedding defined by $\left. \gamma \right|_{I_{\alpha}} = \gamma_{\alpha}$. For each $s \in I$, let $C(s)$ denote the circle in $f_1^{-1}(s)$ containing $\gamma(s)$. Note that $C(s)$ winds around $x$ and, according to the previous paragraph, intersects the complement of each disc $U_i$, $i \neq 1$. Therefore, for each $i \neq 1$, $f_i$ vanishes somewhere on $C(s)$. Since $s \in I$ is a regular value of $f_1$, each $C(s)$ is an integral curve of the Hamiltonian vector field $X_{f_1}$; Denote by $T : I \to (0, \infty)$ the (smooth) function giving the period of the integral curve of $X_{f_1}$ along $C(s)$.

Let's consider the subset of the ``energy-time space''
\[ R := \{ (s,t) \in \R^2 \, : \, s \in I, \, t \in (0, T(s)) \,  \} \, . \]
The map $\Phi : R \to f_1^{-1}(I) \setminus \gamma(I) : (s,t) \mapsto \phi^{f_1}_t(\gamma(s))$ is a diffeomorphism onto its image. We observe that $(T \Phi)(\partial_t) = X_{f_1}$, that $\Phi^*f_1 = s$ and that $-\mathrm{d}s = (\partial_t) \lrcorner \, (\mathrm{d}s \wedge \mathrm{d}t)$, whence $\Phi^*\omega = \mathrm{d}s \wedge \mathrm{d}t$. Moreover, for $s \in I$ and $i \neq 1$, we have
\[ \underset{t \in (0, T(s))}{\mathrm{sup}} \,  (\Phi^*f_i)(s,t) = \underset{C(s)}{\mathrm{max}} \,  f_i \; , \; \; \underset{t \in (0, T(s))}{\mathrm{inf}} \,  (\Phi^*f_i)(s,t) = \underset{C(s)}{\mathrm{min}} \,  f_i = 0 \, ,\]
and $\Phi^*f_i$ oscillates at least twice between those extreme values as $t$ varies in $(0, T(s))$. The rest of the proof is a computation:
\begin{align}
\notag \int_{U_1} \sum_{i=1}^N |\{f_i, f_1\}| \, \omega &= \sum_{i=2}^N \int_{U_1} |\mathrm{d}f_i(X_{f_1})| \, \omega \ge \sum_{i=2}^N \int_{\Phi(R)} |\mathrm{d}f_i(X_{f_1})| \, \omega \\
\notag & = \sum_{i=2}^N \int_{R} |\mathrm{d}(\Phi^*f_i)(\partial_t)| \, \mathrm{d}s \wedge \mathrm{d}t \\
\notag &= \sum_{i=2}^N \int_{s \in I} \int_{t=0}^{T(s)} \left| \frac{d(\Phi^*f_i)}{dt} \right| \, dt ds \\
\notag & \ge \sum_{i=2}^N \int_{s \in I} 2 \, \left( \underset{t \in (0, T(s))}{\mathrm{sup}} \,  \Phi^*f_i - \underset{t \in (0, T(s))}{\mathrm{inf}} \,  \Phi^*f_i \right) \, ds \\
\notag & \ge 2  \int_{s \in I} \sum_{i=2}^N f_i(\gamma(s)) \, ds \\
\notag &= 2 \int_{s \in I} (S_{\F}(\gamma(s)) - f_1(\gamma(s))) \, ds \\
\notag & \ge 2 \int_{s \in I} (1 - s) \, ds = 1 \, .
\end{align}

\QED

\vspace{12pt}

\begin{remark}\label{rmk:oscillation}
We recall that the factor $2$ appearing in the course of the calculation is due to the fact that each map $t \mapsto \Phi^*f_i(s,t)$ goes back and forth between its extreme values at least twice. In a situation where each such map oscillated $2k$ times between its extreme values, the factor $2k$ could be used instead. In \cref{sec:star_sphere}, we shall encounter some functions that oscillate more than twice between their extremes due to the presence of symmetries, thereby allowing us to increase this factor $2$ to the number of those oscillations. This will help explain the lower bound $1$ for $x$ near the set of points at which $\U$ localizes in \Cref{star_ineq}.
\end{remark}\vspace{6pt}

\subsection{Operations on the Poisson bracket function}\label{sec:operations}

Let $\F = \{f_i \}_{i \in I}$ be a locally finite collection of smooth functions on a symplectic manifold $(M, \omega)$, with $I$ a countable index set, we study the behaviour of the functions
\[ S_{\F}(x) = \sum_{i \in I} f_i(x) \; \mbox{ and } \; P_{\F}(x) = \sum_{i, j \in I} |\{f_i, f_j\}| \]
under certain operations.\vspace{12pt}

\emph{Condensation}. Given another countable set $J$ and a surjective map $c : I \to J$, consider the positive collection $\F' = \{f'_j\}_{j \in J}$ obtained by setting 
\[ f'_j := \sum_{i \in c^{-1}(j)} f_i, \quad  \forall \, j \in J \, . \]
Clearly, $S_{\F'} = S_{\F}$. Linearity of the Poisson bracket and the triangle inequality easily imply $P_{\F'} \le P_{\F}$. \vspace{12pt}

\emph{Fragmentation}. Suppose $f_i \in \F$ has disconnected support. Write $\mathrm{supp}(f_i) = A \sqcup B$, where $A$ and $B$ are both unions of connected components of $\mathrm{supp}(f_i)$, and set $f_A := \left. f_i \right|_A$ and $f_B := \left. f_i \right|_B$. Consider the new positive collection $\F' = (\F \setminus \{f_i\}) \cup \{f_A, f_B\}$. Since $f_i = f_A + f_B$ and $\{f_A, f_B\} = 0$, we easily get $S_{\F'} = S_{\F}$ and $P_{\F'} = P_{\F}$. We can evidently iterate this operation on $A$ and $B$; In fact, since each $\mathrm{supp}(f_i)$ has countably many connected components, we can replace $f_i$ in $\F$ by the collection of its restrictions to each of the components of $\mathrm{supp}(f_i)$. Doing so for each $i \in I$, fragmentation leads to a positive collection $\F'$ such that each $f'_i \in \F'$ has connected support, $S_{\F'} = S_{\F}$ and $P_{\F'} = P_{\F}$.\vspace{12pt}

\emph{Lift}. Suppose $p : (M', \omega') \to (M, \omega)$ is a symplectic covering map, \emph{i.e.} a covering map such that $p^*\omega = \omega'$.  Consider the collection $p^*\F := \F' = \{f'_i\}_{i \in I}$ on $M'$ given by $f'_i = p^*f_i$. Note that this collection is locally finite even when $p$ has infinite degree. Clearly, $S_{\F'} = p^*S_{\F}$. Since $p$ is a local symplectic diffeomorphism, we have $p^*\{g,h\} = \{p^*g, p^*h\}$ for all $g,h \in C^{\infty}(M)$; in particular, we deduce $P_{\F'} = p^*P_{\F}$. \vspace{0pt}

\begin{remark}
If $\F$ is subordinate to an open cover $\U$:
\begin{itemize}
 \item A condensation $\F'$ of $\F$ is subordinate to a corresponding union $\U'$ of sets in $\U$, namely $\mathrm{supp}(f'_j) \subset \cup_{i \in c^{-1}(j)} U_j$;
 \item A fragmentation $\F'$ of $\F$ is ``subordinate to $\U$'' in the sense that it is subordinate to an open cover $\U' = \{U'_j\}_{j \in J}$ obtained from $\U$ by counting sets in $\U$ multiple times, namely if $f'_j \in \F'$ comes from fragmenting a certain $f_i \in \F$, then $\mathrm{supp}(f'_j) \subset U'_j := U_i$; 
 \item A lift $p^*\F$ of $\F$ is subordinate to the open cover $p^*\U = \{p^{-1}U_i\}_{i\in I}$ of $M'$.
 \end{itemize}
\end{remark}\vspace{0pt}

\subsection{Star inequality for $M \neq S^2$} \label{sec:star_genus}

\begin{figure}[h!]
\centering
\includegraphics[width=100mm]{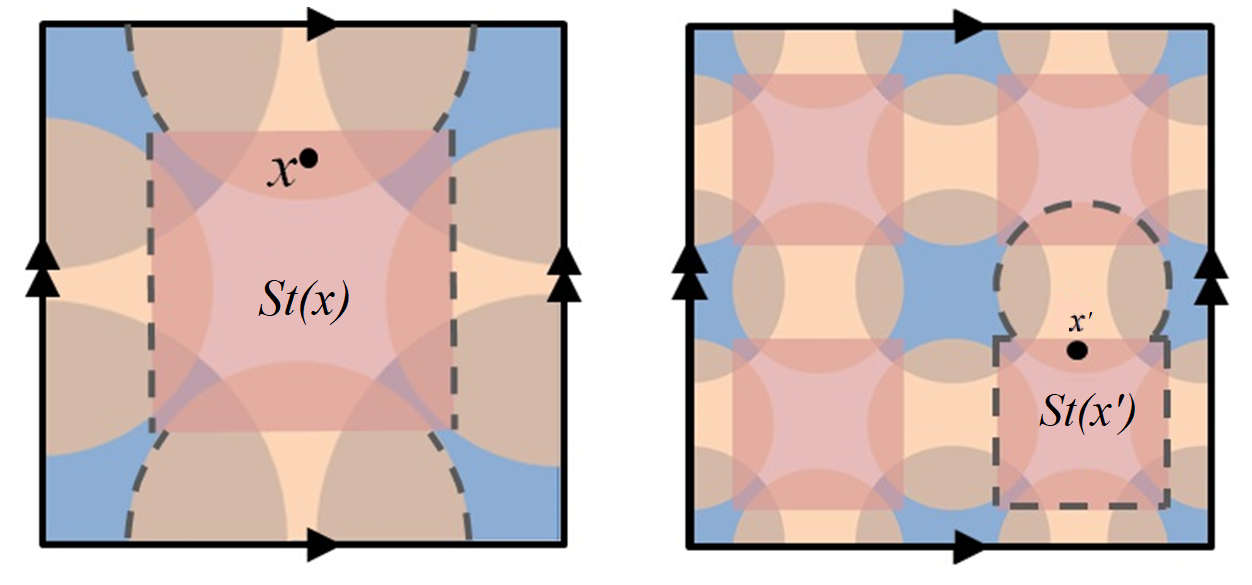}
\caption{Example on the torus of the lift of an open cover along a covering map of degree $4$. Note that $St(x')$ is a confined disc, but not $St(x)$.}\label{fig:star}
\end{figure}

We now prove \Cref{star_ineq} when $M \neq S^2$. We shall use the operations introduced in \cref{sec:operations} to simplify the topology of a star set into that of a confined essential disc (see \Cref{fig:star} for an example of such a simplification). \vspace{12pt}

\vspace{12pt}

\textit{Proof of \Cref{star_ineq}}. Let $p : (M', \omega') \to (M, \omega)$ be a symplectic universal covering map. Since $M \neq S^2$, $M'$ is diffeomorphic to $\R^2$ by the uniformization theorem.

\emph{Step 1 - Lift and fragmentation}. Put $K = \pi_1(M, \ast)$. Since discs are contractible and covering maps have the unique homotopy lifting property, the lift $p^{-1}(U_i)$ of each $U_i \in \U$ along $p$ consists of the disjoint union of discs $U'_{i,k}$ ($k \in K$) each symplectomorphic via $p$ to $U_i$, which we write $p^{-1}(U_i) = \sqcup_{k \in K} U'_{i,k}$. In other words, the open cover $\U'= \{U'_{i,k}\}_{(i,k) \in I \times K} $ is the fragmentation (into connected discs) of the lift of $\U$ to $M'$. Similarly, for each $f_i \in \F$, the function $p^*f_i$ fragments into functions $f'_{i,k}$ respectively supported in $U'_{i,k}$; In other words, the positive collection $\F' = \{f'_{i,k}\}_{(i,k) \in I \times K}$ is a fragmentation of the lift of $\F$ to $M'$, and $\F' \prec \U'$. We note that $\U'$ is still locally finite and in general position.

\emph{Step 2 - Equality of corresponding star integrals}. Let $x \in M$ and pick $x' \in p^{-1}(x)$. Denote simply $St(x) = St(x;\U)$ and $St(x') = St(x';\U')$. We claim that
\[ \int_{St(x)} \sum_{i} \sum_{j : U_j \in \U_x} |\{f_i, f_j\}| \, \omega \;  = \; \int_{St(x')} \sum_{(i,k)} \sum_{(j,l) : U'_{j,l} \in \U'_{x'}} |\{f'_{i,k}, f'_{j,l}\}| \, \omega' \, . \]
Indeed, $U_j \in \U_x$ if and only if there exists (a necessarily unique) $l \in K$ such that $U'_{j,l} \in \U'_{x'}$; We denote $l(j)$ the unique such index. The restriction $\left. p \right|_{U'_{j,l(j)}} : U'_{j,l(j)} \to U_j$ is then a symplectic diffeomorphism. We also note that for $k \neq h$, $\{f'_{i,k}, f'_{j,l}\}$ and $\{f'_{i,h}, f'_{j,l}\}$ have disjoint support, since $f'_{i,k}$ and $f'_{i,h}$ have disjoint support. As a result, we get
\begin{align}
\notag \sum_{(i,k)} & \sum_{(j,l) : U'_{j,l} \in \U'_{x'}}  \int_{St(x')}  |\{f'_{i,k}, f'_{j,l}\}| \, \omega' = \sum_{(i,k)} \sum_{j : U_j \in \U_x} \int_{U'_{j,l(j)}}  |\{f'_{i,k}, f'_{j,l(j)}\}| \, \omega' \\
\notag &= \sum_{i} \sum_{j : U_j \in \U_x} \int_{U'_{j,l(j)}}  \left| \left\{ \sum_{k} f'_{i,k}, f'_{j,l(j)}\right\} \right| \, \omega' \\
\notag &= \sum_{i} \sum_{j : U_j \in \U_x} \int_{U'_{j,l(j)}}  \left| \left\{ p^*f_i, p^*f_j\right\} \right| \, \omega'
\\
\notag &= \sum_{i} \sum_{j : U_j \in \U_x} \int_{U_j}  \left| \left\{ f_i, f_j\right\} \right| \, \omega = \sum_{i} \sum_{j : U_j \in \U_x} \int_{St(x)}  \left| \left\{ f_i, f_j\right\} \right| \, \omega \, .
\end{align}

\textit{Step 3 - Confinement of lifted star}. Fix $x \in M$. For any fixed $x' \in p^{-1}(x)$, the star $St(x')$ is a finite union of open discs in general position. Consequently, $St(x')$ has piecewise smooth boundary, and the boundary components are disjoint contractible piecewise smooth embedded circles. Since $M'$ is diffeomorphic to $\R^2$, the Jordan-Schoenflies theorem implies that each of these boundary components bounds a unique (up to parametrization) embedded open disc, this disc containing any connected component of the complement of its boundary that it intersects. From this, we can argue that there is only one such disc $D'$ that contains $St(x')$. (More details on such argument are given in \Cref{sec:disp_two}.) Clearly, no disc $U' \in \U'_{x'}$ contains $\partial D' \subset \partial St(x')$, and neither does any disc $U' \in \U' \setminus \U'_{x'}$, since this would imply that the curve $\partial D'$ is contractible in $M' \setminus \{x'\}$ and therefore bounds another disc than $D'$, which goes against its uniqueness. We conclude that the component $\partial D'$ of $\partial St(x')$ is contained in no $U' \in \U'$, hence that $St(x')$ is confined.

\textit{Step 4 - Condensation and conclusion}. Let $D''$ be an open disc with smooth boundary which contains $D'$ and which is contained in an arbitrarily small neighbourhood of $D'$, so that no disc $U' \in \U'$ contains $\partial D''$. The cover $\U'$ refines the cover $\U'' := (\U' \setminus \U'_{x'}) \cup \{D''\}$, and $D''$ is confined and essential in the cover $\U''$. The positive collection $\F''$ obtained from $\F'$ by condensation of the functions $\{f'_j : U'_j \in \U'_{x'}\}$, \emph{i.e.} by substitution of these functions by their sum, is subordinate to $\U''$. Hence, by property of condensation, we have
\[ \int_{D''} \sum_{f''_i \in \F''} \sum_{f'_j : U'_j \in \U'_{x'}} |\{f''_i, f'_j\}| \, \omega' \; \ge \; \int_{D''} \sum_{f''_i \in \F''}  \left| \left\{f''_i, \sum_{f'_j : U'_j \in \U'_{x'}} f'_j \right\} \right| \, \omega'  \, , \]
while the right-hand side is at least $1$ by \Cref{conf_ess}. Since the integrand in the left-hand side integral is supported in $St(x') \subset D''$, we conclude from step 2.

\QED \vspace{6pt}

\begin{remark}\label{rmk:oscillation2}
In \cref{sec:star_sphere}, we shall encounter situations where each $U_j \in \U_x$ ($x$ fixed) comes equipped with a $\Z/n\Z$-action, and these actions (indexed by $j$) are compatible in the sense that they coincide on intersections between the discs in $\U$, thereby defining an action on $St(x)$. For each $j$, this induces a $\Z/n\Z$-action on the collection of the restrictions of the functions $f_i$ to $U_j$. All these actions then lift to (compatible) $\Z/n\Z$-actions on the discs $U'_{j,k} \in \U'_{x'}$ and on the restrictions of the functions $f'_{i,l}$ to those discs, and therefore induce a $\Z/n\Z$-action on $St(x')$. Condensation of the discs in $\U_x$ being an equivariant operation with respect to these actions, and taking into account \Cref{rmk:oscillation} in \emph{Step 4}, we deduce in these $\Z/n\Z$-equivariant situations that
\[ \int_{St(x)} \sum_{i} \sum_{j : U_j \in \U_x} |\{f_i, f_j\}| \, \omega \;  \ge n \, . \]
\end{remark} \vspace{6pt}

\vspace{12pt}

\subsection{Star inequality for $M = S^2$} \label{sec:star_sphere}

We now come to the proof of the different statements in \Cref{star_ineq} for $M=S^2$. We shall consider first the case when all star sets in $\U$ are confined, secondly the case when $\U$ is $4$-localized and then the case when $\U$ is $3$-localized. We shall finally prove the statement for $x$ sufficiently close to the points at which $\U$ localizes.  \vspace{12pt}

\textit{Proof of \Cref{star_ineq} when $St(x)$ is confined}. Since $St(x)$ is confined, $\partial \overline{St(x)}$ is nonempty and consists of disjoint piecewise smooth embedded circles, at least one of which is contained in no single $U_i \in \U$. Let $C$ be such a circle. By the Jordan-Schoenflies theorem for $S^2$, $S^2 \setminus C$ consists of two discs. Since $St(x)$ is connected, precisely one of these two discs contains $St(x)$; denote it $D$.

Let $D'$ be an open disc with smooth boundary which contains $D$ and which is contained in an arbitrarily small neighbourhood of $D$, so that no disc $U_i \in \U$ contains $\partial D'$. We form the new cover $\U' = (\U \setminus \U_x) \cup \{D'\}$ and the new positive collection $\F' = (\F \setminus \{f_i : U_i \in \U_x\}) \cup \{\sum_{i : U_i \in \U_x} f_i \}$. Note that $\F' \prec \U'$ and that $D'$ is a confined essential set for $\U'$. By property of condensation, we have
\[ \int_{D'} \sum_{f'_i \in \F'} \sum_{f_j : U_j \in \U_{x}} |\{f'_i, f_j\}| \, \omega \; \ge \; \int_{D'} \sum_{f'_i \in \F'}  \left| \left\{f'_i, \sum_{f_j : U_j \in \U_{x}} f_j \right\} \right| \, \omega  \, , \]
and the right-hand side is at least $1$ by \Cref{conf_ess}. Since the integrand in the left-hand side integral is supported in $St(x) \subset D'$, we get the result.

\QED

\vspace{12pt}

We now turn to the cases when $\U$ is $m$-localized for $m \ge 3$. Analogously to our use of covering maps in the proof of the star inequality in genus $g \ge 1$, we shall make use of ramified covering maps from the torus to the sphere. Let's recall some terminology and facts, loosely following \cite[Chapter 4]{D} and \cite[Chapter 19]{Fu}.  Given two connected closed (oriented) surfaces $M$ and $N$, a map $p : N \to M$ is a \emph{ramified covering map} (or a \emph{branched covering map}) if it is nonconstant and holomorphic for some complex structures on $N$ and $M$ (compatible with the given orientations). Its critical points and critical values are called \emph{ramification points} and \emph{branched points}, respectively. Let $R \subset N$ be the set of ramification points and $B = p(R) \subset M$ be the set of branched points. Both sets are finite and $p : N \setminus p^{-1}(B) \to M \setminus B$ is a $d$-sheeted covering map, where $d$ is the \emph{degree of $p$}. For any $y \in N$, there are complex coordinates $z$ around $y$ and $w$ around $x=p(y)$ with respect to which the map $p$ reads $w = z^e$ for some integer $e = e_p(y) \ge 1$ called the \emph{ramification index of $p$ at the point $y$}. The inequality $e_p(y) > 1$ holds if and only if $y \in R$, and we have $d = \sum_{y \in p^{-1}(x)} e_p(y)$ for all $x \in N$. A ramified covering map of degree $d$ is said to be \emph{simple} if $e_p(y) = d$ for all $y \in R$. The following facts can be deduced from the sort of reasoning given in \cite[Chapter 4.2.2]{D}:  \vspace{12pt}

\textit{Fact A}. Let $D \subset M$ be an open disc containing only one branched point $x$ of $p : N \to M$. Let $D'$ be a connected component of $p^{-1}(D)$. Then $D'$ is an open disc containing only one ramification point $y$ of $p$, and the restriction $p|_{D'} : D' \setminus \{y\} \to D \setminus\{x\}$ is a regular covering map of degree $e_p(y)$. The group $\Z/e_p(y) \Z$ therefore acts on $D'$, fixing $y$ and acting on $D' \setminus \{y \}$ as the deck transformation group of $p|_{D'}$: Its action on each fiber $(p|_{D'})^{-1}(x')$, $x' \in D \setminus \{x\}$, is free and transitive and is generated by the monodromy along a generator of $\pi_1(D \setminus \{x\}, \{x'\}) $.   \vspace{12pt}

\textit{Fact B}. (Riemann's Existence Theorem, \cite[p. 49]{D} ) Let $M$ be a connected closed oriented surface and $\Delta \subset M$ be a finite set. Given $d \ge 1$ and a representation $\rho : \pi_1(M \setminus \Delta) \to S_d$ (the symmetric group on $d$ elements) that acts transitively on a set of $d$ elements, there exist a connected closed oriented surface $N$ and a ramified covering map $p : N \to M$ of degree $d$ with branch locus $B \subset \Delta$ that realises $\rho$ as its monodromy homomorphism. \vspace{12pt}

Let $b, d \ge 1$ be coprime integers. Given a finite set $B \subset S^2$ of cardinality $b+1 \ge 2$, the fundamental group $\pi_1(S^2 \setminus B)$ is isomorphic to the free group on $b$ elements, whose generators $g_1, \dots, g_b$ can be thought of as small embedded circles around $b$ points in $B$. (The element $g_1 \dots g_b$ corresponds to a small embedded circle around the remaining point of $B$.)   Sending $g_1, \dots, g_b$ to the same generator $h$ of $C_d$ (the cyclic group of order $d$) determines a transitive representation $\rho : \pi_1(M \setminus B) \to S_d$ and thus determines a ramified covering map $p : N \to S^2$ of degree $d$. Since $\rho(g_1 \dots g_b) = h^b$ is another generator of $C_d$, it follows that $p$ is simple with branch locus $B$. From the Riemann-Hurwitz formula, the Euler characteristic of $N$ is $b+1 - d(b-1)$, which is nonpositive if and only if $b, d \ge 2$. As special instances of this construction:\vspace{12pt}

\textit{Fact 1}. Given a collection $B$ of $4$ points in $S^2$, there exists a simple ramified covering map $p : T^2 \to S^2$ of degree $2$ with branch locus $B$. \vspace{12pt}

\textit{Fact 2}. Given a collection $B$ of $3$ points in $S^2$, there exists a simple ramified covering map $p : T^2 \to S^2$ of degree $3$ with branch locus $B$. \vspace{12pt}

The crux of the proof of \Cref{star_ineq} lies in the next lemma, which we shall ultimately apply to the covering maps from Facts 1 and 2:

\begin{lemma}\label{lem_star_rami}
Let $(M, \omega)$ be a closed symplectic surface, $\U = \{U_i\}_{i \in I}$ be a finite open cover by discs in general position localized at the points $X = \{x_1, \dots, x_m\}$ and $\F = \{f_i\}_{i \in I} \prec \U$ be a positive collection. Let's assume there exists a simple ramified covering map $p : N \to M$ of degree $d$ with branch locus $B = \{x_1, \dots, x_b\} \subset X$, and let $R \subset N$ be its ramification locus. Then for all $\epsilon > 0$, there exist a finite set $L$, a surjective map $q : L \to I$, a finite open cover $\V = \{V_l\}_{l \in L}$ of $N$ by discs in general position localized at the points of $Y := p^{-1}(X)$, a positive collection $\GG = \{g_l\}_{l \in L} \prec \V$ and a symplectic form $\eta$ on $N$ such that for all $y \in N$,
\[  \int_{M} \sum_{i \in I} \sum_{j : U_j \in \U_{p(y)}} |\{f_i, f_j\}| \, \omega  \, \ge  \, \frac{1}{d} \, \int_{N} \sum_{l \in L} \sum_{h : V_h \in \V_{y}} |\{g_l, g_h\}| \, \eta \, + \, O(\epsilon)  \; .   \]
Moreover, for $y$ sufficiently close to $Y$, we have the equality
\[  \int_{M} \sum_{i \in I} \sum_{j : U_j \in \U_{p(y)}} |\{f_i, f_j\}| \, \omega  \, =  \, C \, \int_{N} \sum_{l \in L} \sum_{h : V_h \in \V_{y}} |\{g_l, g_h\}| \, \eta \, + \, O(\epsilon)  \; ,   \]
where $C=1/d$ if $y$ is close to $R \subset Y$ and $C=1$ if $y$ is close to $Y \setminus R$.
\end{lemma}
\vspace{12pt}

The following technical result will be used in the proof:

\begin{lemma} \label{lem_const_func}
Let $(M, \omega)$ be a symplectic surface, $\U = \{U_i\}_{i \in I}$ be a locally finite open cover by discs in general position and $\F = \{f_i\}_{i \in I} \prec \U$ be a positive collection. Let $x_1, \dots, x_m \in M \setminus \bigcup_{i \in I} \, \partial \overline{U}_i$.  Given any sufficiently small  $\epsilon > 0$, there exist a positive collection $\F' = \{f'_i\}_{i \in I} \prec \U$ and arbitrary small discs $D'_k \supset D_k \ni x_k$ ($1 \le k \le m$) such that each function $f'_i$ is locally constant on $\sqcup_{k=1}^m D_k$ and coincides with $f_i$ outside $\sqcup_{k=1}^m D'_k$. Moreover, $|\{f'_i, f'_j\}| < |\{f_i, f_j\}| + \epsilon$ over $M$ for all $i,j \in I$.
\end{lemma}

\textit{Proof}. For $\delta > 0$ sufficiently small and each $1 \le k \le m$, there is a Darboux chart $\Phi_k : (D'_k, \omega) \to (B^2(\delta), \omega_0)$ sending $x_k \in D'_k \subset M$ to $0 \in B^2(\delta) \subset \R^2$, where $\omega_0$ denotes the standard symplectic form and $B^2(\delta)$ denotes the round closed ball of radius $\delta$ centred at $0$. We shall take $\delta$ so small that the different discs $D'_k$ are pairwise disjoint, are included inside any $U_i$ they intersect, and are such that the following inequality holds (which is possible since each function $P_{ij} := |\{f_i, f_j\}|$ is continuous)
\[ \underset{1 \le k \le m}{\mathrm{max}} \, \underset{i,j \in I}{\mathrm{max}} \, \underset{x,y \in D'_k}{\mathrm{max}} |P_{ij}(x) - P_{ij}(y)| < \frac{\epsilon}{2} \, . \]

Let $C = 1+ \mathrm{max}_k \, \mathrm{max}_{i,j \in I} \, \mathrm{max}_{x \in D'_k} \, P_{ij}(x)$ and take $\epsilon < C$. Let $\sigma > 0$ be smaller than $\epsilon \delta/8C < \delta/4$. Fix also a smooth nondecreasing function  $\rho : [0, \delta] \to [0, \delta]$ such that $\rho(r) \le r$ over $[0, \delta]$, $\rho(r) = 0$ for $r \in [0, \sigma]$, $\rho(r) = r$ for $r \in [\delta - \sigma, \delta]$ and $\rho'(r) < 1 + 4\sigma/\delta < 1 + \epsilon/2C$ over $[0, \delta]$.  Using polar coordinates $(r, \theta)$ on $B^2(\delta)$, consider
\[\Psi : B^2(\delta) \to B^2(\delta) : (r, \theta) \mapsto (\rho(r), \theta) \, .\]
Given a smooth function $f : B^2(\delta) \to \R$, the pullback $\Psi^*f : B^2(\delta) \to \R$ is smooth, is constant on $B^2(\sigma)$ and equals $f$ near $\partial B^2(\delta)$. Given two functions $f, g : B^2(\delta) \to \R$, a straightforward calculation yields $|\{\Psi^* f, \Psi^* g\}_0| = |\rho'(r)| \, \left| \frac{\rho(r)}{r} \right| \, |\Psi^*\{f,g\}_0| <  (1 + \epsilon/2C) \, |\Psi^*\{f, g\}_0|$, where all functions are evaluated at $(r, \theta)$ and $\{ \cdot , \cdot \}_0$ denotes the Poisson bracket with respect to $\omega_0 = r dr \wedge d\theta$.

We now define a smooth map $\psi : M \to M$ as equal to the identity outside $\sqcup_{k=1}^m D'_k$ and as equal to $\psi|_{D'_k} = \Phi_k^{-1} \circ \Psi \circ \Phi_k : D'_k \to D'_k$ for all $1 \le k \le m$. For each $i \in I$, we define the smooth function $f'_i = \psi^*f_i$. Setting $D_k := \Phi_k^{-1} B^2(\sigma)$, each $f'_i$ is seen to be locally constant on $\sqcup_{k=1}^m D_k$. The collection $\F'$ is still positive and subordinate to $\U$. It also follows that for all $i,j \in I$, $|\{f'_i, f'_j\}| = |\{f_i, f_j\}|$ outside $\sqcup_{k=1}^m D_k$, whereas within this set, we have
\[|\{f'_i, f'_j\}| < \left(1 + \frac{\epsilon}{2C} \right) \, |\psi^*\{f_i, f_j\}| \le |\psi^*\{f_i, f_j\}| + \frac{\epsilon}{2} \le |\{f_i, f_j\}| + \epsilon \, .\]

\QED

\vspace{12pt}

\textit{Proof of \Cref{lem_star_rami}}. We split the proof into several steps.

\textit{Step 1 - Perturbation of $\F$}. By \Cref{lem_const_func}, for any sufficiently small $\epsilon > 0$, we may find small discs $x_k \in D_k $ ($1 \le k \le b$) and a positive collection $\F' \prec \U$  whose functions $f'_i$ are constant on each $D_k$ and which satisfies $|\{f'_i, f'_j\}| < |\{f_i, f_j\}| + \epsilon$ for all $i,j \in I$. We shall assume that each $D_k$ is small enough to intersect only the discs $U_i \in \U_{x_k}$ and that it is contained in all of those.

\textit{Step 2 - Choice of $L$ and $q$}. Since the branch locus $B$ of $p$ is assumed to be a subset a the localization locus $X$ of $\U$, all $U_i \in \U$ intersect $B$ in at most one point. Let's write $I = I_0 \sqcup I_1$, where $i \in I_0$ if and only if $U_i \cap B = \emptyset$. For $h = 1, \dots, d$, we set $L_{0,h}  := \{(i, h)\}_{i \in I_0} \cong I_0$, $L_0 := \sqcup_{h=1}^d L_{0,h}$, $L_1 :=  \{(i, R)\}_{i \in I_1} \cong I_1$ and $L :=  L_0 \sqcup L_{1}$. In other words, $L$ is the union of $d$ copies of $I_0$ and one copy of $I_1$, and the second index $h=1,\dots, d, R$ in the elements $(i,h) \in L$ helps keeping track of which copy of which $I_j$ the index $i$ belongs to. We define $q : L \to I : (i,h) \mapsto i$. When there is little risk for confusion, we may drop the index $R$ from the notation.

\textit{Step 3 - Choice of $\V$}. We define the cover $\V = \{V_l\}_{l \in L}$ on $N$ as given by all the connected components of the sets $p^{-1}(U_i)$, for $i \in I$. More precisely, for $i \in I_0$, $p^{-1}(U_i)$ is the disjoint union of $d$ discs $V_{i,h}$ ($h = 1 \dots, d$), all diffeomorphic to $U_i$ under $p$ and containing no point of $R$. For $i \in I_1$, the intersection $U_i \cap B$ is nonempty, say equal to $\{x_k\}$ for some $1 \le k \le b$; Then \emph{Fact A} implies that $p^{-1}(U_i)$ is a single disc $V_{i,R}$ and that $p|_{V_{i,R}} : V_{i,R} \to U_i$ is a ramified covering map of degree $d$ with a single ramification point $\{y_k\} = p^{-1}(x_k)$.

The cover $\V$ is clearly in general position. Moreover,  $\V$ is localized at the points of $Y := p^{-1}(X) = \bigcup_{k=1}^b \{y_k\} \cup \bigcup_{k>b}^m \bigcup_{h=1}^d \{y_{k,h}\}$, where for $b < k \le m$, $ \{y_{k,1}, \dots, y_{k,d}\} := p^{-1}(x_k)$ with $y_{k,h} \in V_{i,h}$ if $x_k \in U_i$.

We also remark that for each $1 \le k \le b$, the set $p^{-1}(D_k)$ is a disc $D'_k$ that intersects only those $V_l$ that contains $y_k$. $D'_k$ is therefore contained in all those $V_l$.

\textit{Step 4 - Choice of $\GG$}. We define the positive collection $\GG = \{g_l\}_{l \in L} \prec \V$ as the fragmentation of the functions $p^*f'_i$. More precisely, for $i \in I_0$, we have $p^*f'_i =\sum_{h=1}^d  g_{i,h}$ where $ g_{i,h} = (p^*f'_i)|_{V_{i,h}}$, whereas for $i \in I_1$, we have $g_{i,R} = p^*f'_i$. All functions $g_l \in \GG$ are constant on each $D'_k$.

\textit{Step 5 - Choice of $\eta$}. The closed form $p^* \omega$ fails to be symplectic precisely on the discrete ramification locus $R$. It is clear that there is a symplectic form $\eta$ that differs from $p^*\omega$ only inside $\sqcup_{k=1}^b D'_k$. (Although any such form suffices for our purpose, there is some appeal to take $\eta$ such that $\int_{D'_k} \eta = \int_{D'_k} p^* \omega = d \int_{D_k} \omega$. Since each $D'_k$ is contained in every $V_l$ it intersects, it then follows that the $\eta$-area of every $V_l$ coincides with its $p^*\omega$-area. Namely, for $l \in L$, $\int_{V_l} \eta = d \int_{U_{q(l)}} \omega$ if $l \in L_{1}$, whereas $\int_{V_l} \eta =  \int_{U_{q(l)}} \omega$ otherwise.  Hence $c(\V) = d \, c(\U)$ and $\mathrm{Area}(N, \eta) = d \, \mathrm{Area}(M, \omega)$.)

For all $l,h \in L$, we have $\{g_l, g_h\}_{\eta}(y) = \{f'_{(q(l)}, f'_{q(h)}\}_{\omega}(p(y))$. Indeed, both sides of the equality vanish for $y \in D'_k$ and $p(y) \in D_k$, since the functions $g_l$ and $f'_i$ are constant on $D'_k$ and $D_k$ respectively, whereas for $y \not \in \sqcup_{k=1}^b D'_k$, we have $\eta = p^*\omega$ and $p^*f'_{q(l)}(y) = g_l(y)$ for $y \in \mathrm{supp} \, g_l$. Consequently, $|\{p^*f'_i, p^*f'_j\}_{\eta}| = p^*|\{f'_i, f'_j\}_{\omega}|$ over the whole of $N$.

\textit{Step 6 - Star inequality for general $y$}. Fix $y \in N$. We note that  there is a bijection between $\V_y := \{ V_l \in \V \, : \, y \in V_l \}$ and $\U_{p(y)} =  \{ U_i \in \U \, : \, p(y) \in U_i \}$. We denote the bijection $s : \{ i \in I \, : \, U_i \ni p(y) \} \to \{ l \in L \, : \, V_l \ni y \}$, so that $q(s(i)) = i$. We have
\begin{align}
\notag  \int_{N} \sum_{l \in L} \sum_{h : V_h \in \V_{y}} |\{g_l, g_h\}| \, \eta &= \int_{N} \sum_{l \in L} \sum_{j : U_j \in \U_{p(y)}} |\{g_l, g_{s(j)}\}| \, \eta \\
\notag &= \sum_{j : U_j \in \U_{p(y)}}  \int_{V_{s(j)}}  \sum_{i \in I} |\{p^*f'_i, p^*f'_j\}| \, \eta \, ,
\end{align}
where to obtain the second equality, we used that (i) for $i \in I_0$, $p^*f'_i = \sum_{h=1}^d g_{i,h}$  with the functions $g_{i,1}, \dots, g_{i,d}$ supported in disjoint sets, so as to reformulate the sum over $l \in L$ as a sum over $i \in I$, and (ii) each $g_{s(j)}$ is supported in $V_{s(j)}$ and equals $p^*f'_j$ on this set. Now, each $p|_{V_{s(j)}} : V_{s(j)} \to U_j$ is a (possibly ramified) covering map of degree $d(j) = 1$ or $d(j) = d$ depending on whether $j \in I_0$ or $j \in I_1$, respectively. Hence
\[   \int_{V_{s(j)}}  \sum_{i \in I} |\{p^*f'_i, p^*f'_j\}| \, \eta = d(j) \, \int_{U_j} \sum_{i \in I} |\{f'_i, f'_j\}| \, \omega \, . \]
Using that $d(j) \le d$ for all $j$, we obtain
\begin{align}
\notag  \int_{N} \sum_{l \in L} \sum_{h : V_h \in \V_{y}} |\{g_l, g_h\}| \, \eta &\le d \,  \sum_{j : U_j \in \U_{p(y)}}   \, \int_{U_j} \sum_{i \in I} |\{f'_i, f'_j\}| \, \omega \\
\notag &= d\, \int_{M} \sum_{i \in I}  \sum_{j : U_j \in \U_{p(y)}}  |\{f'_i, f'_j\}| \, \omega \\ 
\notag &\le d \, \int_{M} \sum_{i \in I}  \sum_{j : U_j \in \U_{p(y)}}  |\{f_i, f_j\}| \, \omega \, + \, O(\epsilon) \, ,
\end{align}
proving the claimed inequality.

\textit{Step 7 - Star inequality for $y$ near $Y$}. If $y$ is sufficiently close to $Y$, say to $y' \in Y$, then $\V_y = \V_{y'}$ and $\U_{p(y)} = \U_{p(y')}$, hence $St(y) = St(y')$ and $St(p(y)) = St(p(y'))$. Hence it suffices to prove the result for $y \in Y$. Now, because $\V$ is localized at $Y$ and $\U$ is localized at $X = p(Y)$, and because the ramification and branch loci satisfy $R \subset Y$ and $B \subset X$, it follows that $d(j) = e_p(y)$ over the set $\{ j : U_j \in \U_{p(y)} \}$, that is, $d(j)$ equals $1$ if $y \in Y \setminus R$ and equals $d$ if $y \in R$.

For $y \in Y \setminus R$, the computation from the previous step leads to
\begin{align}
 \notag  \int_{N} \sum_{l \in L} \sum_{h : V_h \in \V_{y}} |\{g_l, g_h\}| \, \eta &= \int_{M} \sum_{i \in I}  \sum_{j : U_j \in \U_{p(y)}}  |\{f'_i, f'_j\}| \, \omega \\
\notag &=  \int_{M} \sum_{i \in I}  \sum_{j : U_j \in \U_{p(y)}}  |\{f_i, f_j\}| \, \omega \, + \, O(\epsilon) \, .
\end{align}
For $y \in R$, the same computation rather leads to
\[ \notag  \int_{N} \sum_{l \in L} \sum_{h : V_h \in \V_{y}} |\{g_l, g_h\}| \, \eta = d \, \int_{M} \sum_{i \in I}  \sum_{j : U_j \in \U_{p(y)}}  |\{f_i, f_j\}| \, \omega \, + \, O(\epsilon) \, . \]

\QED

\vspace{12pt}

\textit{Proof of \Cref{star_ineq} for $4$-localized covers on $S^2$}. Let $\U = \{U_i\}_{i \in I}$ be a finite open cover on $S^2$ by discs in general position localized at the points $X = \{x_1, \dots, x_m\}$ with $m \ge 4$, and $\F = \{f_i\}_{i \in I} \prec \U$. By \emph{Fact 1}, there exists a simple ramified covering $p : T^2 \to S^2$ of degree $2$ with branch locus $B = \{x_1, x_2, x_3, x_4\}$. Given $\epsilon > 0$, let $\V = \{V_l\}_{l \in L}$ be the open cover, $\GG = \{g_l\}_{l \in L} \prec \V$ be the positive collection and $\eta$ be the symplectic form on $T^2$ given by \Cref{lem_star_rami}. Then, for any $x \in S^2$, writing $x = p(y)$, we have
\begin{align}
\notag  \int_{St(x)} \sum_{i \in I} \sum_{j : U_j \in \U_{x}} |\{f_i, f_j\}| \, \omega  \, & \ge  \, \frac{1}{2} \, \int_{St(y)} \sum_{l \in L} \sum_{h : V_h \in \V_{y}} |\{g_l, g_h\}| \, \eta \, + \, O(\epsilon)  \\
\notag &\ge \frac{1}{2} + O(\epsilon) \underset{\epsilon \to 0}{\longrightarrow} \frac{1}{2} \, .
\end{align}
Indeed, the first inequality follows from the inequality given in \Cref{lem_star_rami}, while the second inequality follows from \Cref{star_ineq} applied to $(T^2, \eta)$ and the data $\V$ and $\GG$. This proves that \cref{eq_star_ineq2} holds with $C=1/2$.

It remains to prove that \cref{eq_star_ineq2} holds with $C=1$ when $x$ is close to $X$. Let $R \subset Y \subset T^2$ denote the ramification locus of $p$. If $x = p(y)$ is close to $X \setminus B$, that is if $y$ is close to $Y \setminus R$, then
\begin{align}
\notag  \int_{St(x)} \sum_{i \in I} \sum_{j : U_j \in \U_{x}} |\{f_i, f_j\}| \, \omega  \, & =  \,  \int_{St(y)} \sum_{l \in L} \sum_{h : V_h \in \V_{y}} |\{g_l, g_h\}| \, \eta \, + \, O(\epsilon)  \\
\notag &\ge 1 + O(\epsilon) \underset{\epsilon \to 0}{\longrightarrow} 1 \, .
\end{align}
This follows from the same argument as before, except that the first equation now results from the appropriate equality in \Cref{lem_star_rami}.

For $x = p(y)$ close to $B \subset X$, that is for $y$ close to $R \subset Y$, the equality in \Cref{lem_star_rami} yields
\begin{align}
\notag  \int_{St(x)} \sum_{i \in I} \sum_{j : U_j \in \U_{x}} |\{f_i, f_j\}| \, \omega  \, & =  \, \frac{1}{2} \, \int_{St(y)} \sum_{l \in L} \sum_{h : V_h \in \V_{y}} |\{g_l, g_h\}| \, \eta \, + \, O(\epsilon)  \, .
\end{align}
We claim that $\int_{St(y)} \sum_{l \in L} \sum_{h : V_h \in \V_{y}} |\{g_l, g_h\}| \, \eta \, \ge \, 2$. Indeed, from \emph{Fact A}, there is an action of $\Z/e_p(y)\Z$ on each disc $V_h \in \V_{y}$ with $e_p(y) = d$ (here, $d =  2$), and these different actions coincide on the intersection of different discs (since they result from the same monodromy action), thereby inducing a $\Z/d\Z$-action on $St(y)$. This action leaves $p^*\omega$ invariant on this set. It also induces an action on the restrictions of the functions $g_l$ to $St(y)$, fixing the functions $g_{i,R}$ with $i \in I_1$ and (cyclically) permuting the functions $g_{i,1}, \dots, g_{i,d}$ for $i \in I_0$ (we are here using notations introduced during the proof of \Cref{lem_star_rami}).

What we need at this point are extensions of the proofs of \Cref{star_ineq} for $T^2$ and of \Cref{conf_ess} to this equivariant setting. To avoid overburdening the presentation with further notation, let's simply consider the case when $St(y)$ is already contained in a disc confined with respect to $\V$, in the sense that its boundary is contained in no $V_l$.  (See however \Cref{rmk:oscillation2} for information about what is involved in the general argument, which would involve lifting everything to the universal cover of $T^2$.) Setting $g_y :=  \sum_{h : V_h \in \V_{y}} g_h$, $g'_i = \sum_{h = 1}^d g_{i,h}$ for $i \in I_0$ and $g'_i = g_{i,R}$ for $i \in I_1$,  we observe that these are all $\Z/d\Z$-invariant and that
\[ \sum_{l \in L} \sum_{h : V_h \in \V_{y}} |\{g_l, g_h\}| \, \ge \, \sum_{i \in I} |\{ g'_{i} , g_y\}|  \, . \]
Recalling the proof of \Cref{conf_ess}, we deduce that around appropriate level-sets of $g_y$, the functions $g'_{i}$ oscillate at least $2d$-times between their extremal values, hence (taking into account \Cref{rmk:oscillation})
\begin{align}
\notag \int_{St(y)} \sum_{l \in L} \sum_{h : V_h \in \V_{y}} |\{g_l, g_h\}| \, \eta \, &\ge \, \int_{St(y)} \sum_{i \in I} |\{g'_i, g_h\}|  \, \eta \, \ge \, d \, ,
\end{align}
with $d = 2$ in the present proof. Consequently,
\begin{align}
\notag  \int_{St(x)} \sum_{i \in I} \sum_{j : U_j \in \U_{x}} |\{f_i, f_j\}| \, \omega  \, & \ge  \, \frac{2}{2} + \, O(\epsilon) \underset{\epsilon \to 0}{\longrightarrow} 1 \, .
\end{align}
This proves that \cref{eq_star_ineq2} holds with $C=1$ when $x$ is close to $X$.

\QED

\vspace{12pt}

\textit{Proof of \Cref{star_ineq} for $3$-localized covers on $S^2$}.  Let $\U = \{U_i\}_{i \in I}$ be a finite open cover on $S^2$ by discs in general position localized at the points $X = \{x_1, x_2, x_3\}$, and $\F = \{f_i\}_{i \in I} \prec \U$. By \emph{Fact 2}, there exists a simple ramified covering $p : T^2 \to S^2$ of degree $3$ with branch locus $B =X$. Arguing along the same lines as in the case of 4-localized covers, but this time with $d=3$, we obtain that for any $x \in S^2$,
\begin{align}
\notag  \int_{St(x)} \sum_{i \in I} \sum_{j : U_j \in \U_{x}} |\{f_i, f_j\}| \, \omega  \, & \ge  \,  \frac{1}{3} \, ,
\end{align}
whereas for $x$ close to $X$,
\begin{align}
\notag  \int_{St(x)} \sum_{i \in I} \sum_{j : U_j \in \U_{x}} |\{f_i, f_j\}| \, \omega  \, & \ge  \, \frac{3}{3} = 1 \, .
\end{align}

\QED

\vspace{12pt}

%%%%%%%%%%%%%%%%%%%%%%%%%%%%%%%%%%%%5

\section{Sharpness of the results} \label{sec:sharpness}

\subsection{Proof of \Cref{Linfty_Lone}}\label{sec:norms}

We shall need the following lemma, which morally results from applying Moser's argument to a path of symplectic forms between $\omega$ and (a rescaling of) $P_{\F} \, \omega$.\vspace{0pt}

\begin{lemma}\label{lower}
Let $(M, \omega)$ be a closed surface, $\U = \{U_i\}_{i \in I}$ a finite open cover on $M$ by discs in general position and $\F = \{f_i\}_{i \in I} \prec \U$ a positive collection. For any $\eta > 0$, there is a diffeomorphism $\phi \in \mathrm{Diff}_0(M)$ such that the positive collection $\F' := \phi^{\ast}\F = \{\phi^*f_i\}_{i \in I}$ subordinate to the cover $\U' := \phi^{\ast}\U = \{\phi^{-1}U_i\}_{i \in I}$  together satisfy
\[ \| P_{\F'} \| \mathrm{Area}(M, \omega) \in \int_M P_{\F}\: \omega \, + \,(-\eta , \eta) \, \]
and
\[\| P_{\F'} \|c(\U') \in \underset{i \in I }{\mathrm{max}} \, \int_{U_i} P_{\F} \: \omega \; + (-\eta , \eta) \; .  \]
\end{lemma}\vspace{6pt}

\textit{Proof}. If $P_{\F} \equiv 0$, we may take $\phi = Id$; We henceforth assume $\int_M P_{\F} \, \omega > 0$. Consider $0 < \eta < E := 2(\int_M \omega \, + \, \int_M P_{\F} \: \omega / \int_M \omega)$ and set $\epsilon = \eta/E$. For $\delta > 0$, by the Stone-Weierstrass theorem, there is a smooth function $P : M \to \R$ satisfying $P_{\F} < P < P_{\F} + \delta$. In particular, $P$ is strictly positive everywhere on $M$. We shall pick  $\delta = \epsilon \cdot \mathrm{min} \left\{ 1 \, , \, \int_M P_{\F} \: \omega/(\int_M \omega)^2\right\}$. Let's consider the differential forms 
\[ \omega_0 = \frac{\int_M \omega}{\int_M P \: \omega}\, P \; \omega \; \mbox{ and } \; \omega_1 = \omega \; .\]

\noindent The two-form $\omega_0$ is well-defined and nondegenerate as $P > 0$, and it is closed as we work on a surface. We observe that $\omega_0$ and $\omega_1$ give the same area to $M$, so that $[\omega_0] = [\omega_1] \in H^2_{dR}(M;\R)$. For $t \in [0,1]$ we set $\omega_t = (1-t) \omega_0 + t \omega_1$, which is path of symplectic forms in the same cohomology class. Moser's path argument then yields a $\phi \in \mathrm{Diff}_0(M)$ such that $\phi^{\ast}\omega_0 = \omega_1$. As a result, for any $\alpha, \beta \in C^{\infty}(M; \R)$,
\begin{align}
\notag \left| \left\{  \phi^{\ast}\alpha, \phi^{\ast} \beta  \right\}_{\omega} (p)  \right| &= \left| \left\{  \alpha,  \beta  \right\}_{\omega_0} (\phi(p))  \right|  = \frac{\int_M P \: \omega}{\int_M \omega} \; \frac{\left| \{ \alpha, \beta \}_{\omega}(\phi(p)) \right|}{P(\phi(p))} \; .
\end{align}
Taking $\alpha = f_i$ and $\beta = f_j$ and summing over all $i$ and $j$,  this yields
\[ P_{\F'}(p) =  \frac{\int_M P \: \omega}{\int_M \omega} \; \frac{P_{\F}(\phi(p))}{P(\phi(p))} \, .\]
Since $\delta \le \epsilon$, we have $P_{\F'} < \int_M (P_{\F} +\delta) \: \omega / \int_M \omega \le \int_M P_{\F} \: \omega / \int_M \omega + \epsilon$ and $\| P_{\F'} \| > ( \int_M P_{\F} \: \omega / \int_M \omega) (1 - \delta/pb(\F)) \ge \int_M P_{\F} \: \omega / \int_M \omega - \epsilon$. Hence
\[ \| P_{\F'} \| \in \frac{\int_M P_{\F} \:  \omega}{\int_M \omega} \, + (-\epsilon, \epsilon) \; \, ,\]
which readily implies the first claimed inequality. We also have
\begin{align}
\notag \mathrm{Area}(\phi^{-1}U_i, \omega) &= \int_{\phi^{-1}U_i} \omega =  \int_{U_i} \omega_0 = \frac{\int_{U_i} P \; \omega}{\int_M P \: \omega}\int_M \omega \, ,
\end{align}
which implies, since $\delta < \int_M P_{\F} \: \omega/(\int_M \omega)^2$,
\begin{align}
\notag \mathrm{Area}(\phi^{-1}U_i, \omega) &< \frac{\int_{U_i} P_{\F} \; \omega}{\int_M P_{\F} \: \omega} \, \int_M \omega \,  + \, \delta \, \frac{(\int_M \omega)^2}{\int_M P_{\F} \: \omega} \le  \frac{\int_{U_i} P_{\F} \; \omega}{\int_M P_{\F} \: \omega} \, \int_M \omega \,  + \, \epsilon
\end{align}
and
\begin{align}
\notag c(\U') &  >  \frac{\underset{i \in I}{\mathrm{max}} \, \int_{U_i} P_{\F} \; \omega}{\int_M P_{\F} \: \omega} \, \int_M \omega \, - \, \delta \, \frac{(\int_M \omega)^2}{\int_M P_{\F} \: \omega} \ge  \frac{\underset{i \in I}{\mathrm{max}} \, \int_{U_i} P_{\F} \; \omega}{\int_M P_{\F} \: \omega} \, \int_M \omega \,  - \, \epsilon \, .
\end{align}
Therefore,
\[   c(\U') \in  \frac{\underset{i \in I}{\mathrm{max}} \, \int_{U_i} P_{\F} \; \omega}{\int_M P_{\F} \: \omega} \, \int_M \omega \,  + \, (-\epsilon, \epsilon)  \, . \]
Multiplying our estimates for $\| P_{\F'} \|$ and $c(\U')$ yields the second claimed inequality.

\QED \vspace{12pt}

\textit{Proof of \Cref{Linfty_Lone}}. The implication `(2) $\Rightarrow$ (1)' is obvious. The implication `(1) $\Rightarrow$ (2)' follows from \Cref{lower}. Indeed, assume (1) holds for some constants $C, C' > 0$. Fix an open cover $\U = \{U_i\}_{i \in I} \in \mathfrak{U}$ and a positive collection $\F \prec \U$. The lemma states that for each $\eta > 0$, there exists $\phi \in \mathrm{Diff}_0(M)$ such that the open cover $\U' = \phi^*\U \in \mathfrak{U}$ and the positive collection $\F' = \phi^*\F \prec \U'$ satisfy
\[ \int_M P_{\F} \, \omega \, \ge \, pb(\F') \, \mathrm{Area}(M, \omega) - \eta \, \ge \, C - \eta \]
and
\[ \underset{i \in I}{\mathrm{max}} \, \int_{U_i} P_{\F} \, \omega \, \ge \, pb(\F') \, c(\U') - \eta \, \ge \, C' - \eta \, . \]
As this is true for all $\eta > 0$, we obtain the statement (2), with the same constants $C$ and $C'$.

\QED \vspace{0pt}

\begin{remark}\label{rmk:fixedpoint}
The diffeomorphism $\phi$ in \Cref{lower} is far from unique, since postcomposing it with a symplectic diffeomorphism leads to another appropriate diffeomorphism. Since $\mathrm{Symp}_0(M, \omega)$ acts transitively on finite sets in $M$, we can require for the diffeomorphism $\phi$ to fix some number of points in $M$. Accordingly, \Cref{Linfty_Lone} also holds for collections $\mathfrak{U}$ that are invariant under the action of the stabilizer subgroup $G \subset \mathrm{Diff}_0(M)$ of any given finite subset of $M$.
\end{remark}\vspace{0pt}

\subsection{Sharpness of the ``confined essential disc'' inequality}\label{sec:sharp_conf_ess}

In view of the proof of \Cref{conf_ess}, it seems feasible to come up with an example of a positive collection for which the inequalities encountered in the course of the proof are nearly equalities. One inequality however appears more difficult to turn into an equality than the others: In the notations of the proof, we would need to find a curve $\gamma$ such that each function $f_i$, $i \neq 1$, attains its maximum on $C(s)$ at $\gamma(s)$, and that for every $s \in I$. Notice however that the choice of $\gamma(s) \in C(s)$ was arbitrary, hence averaging over this choice of point for each $s \in I$, we see that we would need to have, for most $s \in I$ and each $i \neq 1$,
\[  \underset{C(s)}{\mathrm{max}} \,  f_i -  \underset{C(s)}{\mathrm{min}} \,  f_i =  \underset{C(s)}{\mathrm{max}} \,  f_i \; \approx \; \frac{1}{T(s)} \int_0^{T(s)} (\Phi^*f_i)(t,s) \, dt \; .\]
In other words, all functions $f_i$ should be approximately equal to their maximum over most of the interval $t \in [0, T(s)]$, and that for most $s \in [0,1]$. This observation inspires the following proof.\vspace{12pt}

\textit{Proof of  \Cref{prop:sharp_conf_ess_disc}}. Consider the round sphere $S^2 = \{(x,y,z) \in \R^3 \, : \, x^2+y^2+z^2 = 1 \}$ equipped with the usual area form $\omega$. Let $d \ge 2$ be an integer, $h : [-1,1] \to [0,1]$ be a smooth nondecreasing function with $h(u) = 0$ near $u=0$ (and thus for $u < 0$) and $h(u) = 1$ near $u=1$, and $w : \R/2\pi\Z \to [0,1] $ be a smooth function such that $w(t) = 0$ near $t=0$, $w(t) = 1$ for $|t| > 2\pi/3(d+1)$ and monotone over the two remaining intervals (here, $t \in [- \pi, \pi)$).  Consider the collection $\F = \{f_+, f_-, f_0, \dots, f_d\}$ given in cylindrical coordinates $(\theta, z) \in \R/2\pi\Z \times [-1,1]$ as follows:
\begin{align}
\notag &f_+(\theta,z) = h(z), \, f_-(\theta, z) = h(-z) \mbox{ and } \\
\notag & f_j(\theta, z) = \frac{1}{d} \, (1-h(|z|))w\left(\theta + \frac{2\pi j}{d+1}\right) \; \forall \, j =0, \dots, d \, .  
\end{align}
The open cover $\U$ is formed by discs which are slight enlargements of the support of these functions. The north pole is covered by a unique disc, which is therefore essential, and it is also confined.  Using that  $w(\theta + 2\pi j/(d+1))$ equals $1$ on the support of the derivative of $w(\theta + 2\pi k/(d+1))$ whenever $j \neq k$, we get $S_F(\theta, z) \in [1, 1 + 1/d]$, so $\F$ is positive. A straightforward computation yields
\[ \int_M \sum_{f \in \F} |\{f_+, f\}| \, \omega = \int_M \sum_{j=0}^d |\{f_+, f_j\}| \, \omega = 1 + \frac{1}{d} \, .  \]
Letting $d \to \infty$ proves \Cref{prop:sharp_conf_ess_disc} in the case of the sphere.

The closed disc $D$ given by the support of the function $f_+$ can be embedded into any surface $M$ as a disc $D'$. For each $d \ge 2$, the restrictions of the functions $f_+$ and $f_j$ ($j = 0, \dots, d$) to $D$, interpreted as defined on $D'$, can be smoothly extended to $M$ to be supported in a small neighbourhood of $D'$. By adding functions supported in discs contained in the complement of $D'$, we can form a positive collection $\F$ on $M$ subordinate to an open cover $\U$ by discs such that $U_1 := D'$ is confined and essential.  This concludes the proof of  \Cref{prop:sharp_conf_ess_disc}.

\QED

\vspace{12pt}

\begin{remark}
Back to the case of the sphere, a slightly more involved calculation yields $\int_M \, P_{\F} \, \omega = 8(1 + 1/d)$, independently of the details of the functions $h$ and $w$. This is bigger than the lower bound $2$ guaranteed by \Cref{cor_conf_ess}.  By choosing $h$ appropriately, the minimal value for $c(\U)$ we can arrange to have is approximately $\mathrm{Area}(M, \omega)/3$, so that $\int_M \, P_{\F} \, \omega \, \gtrsim (8/3)(\mathrm{Area}(M, \omega)/c(\U))$. The constant $8/3$ is certainly bigger than the constant obtained $1/2$ by \cite{BLT} and than the best constant we obtained for covers of surfaces of higher genus, namely $1$.
\end{remark}

\subsection{Improving the lower bound in the Poisson bracket theorem}\label{sec:sharp_pb_thm}

We now prove \Cref{prop:lowerbound2} which, by improving the lower bound on $\int_U P_{\F} \, \omega$ when $U \in \U$ is a confined essential disc, suggests that the lower bounds in \Cref{pb_thm} might fail to be sharp for all open covers. \vspace{12pt}

\textit{Proof of \Cref{prop:lowerbound2}}. Let $U_k \in \mathcal{J}_c(\U)$. Since $\{f, 1\} = 0$ for all $f$, we have
\begin{align}
\notag \sum_{i,j \in I} |\{f_i, f_j\}| &= \sum_{i \in I} |\{f_i, f_k\}| + \sum_{i \in I} \sum_{j \in I \setminus \{k\}} |\{f_i, f_j\}| \\
\notag &\ge  \sum_{i \in I} |\{f_i, f_k\}| + \sum_{i \in I} \left| \left\{ f_i, \sum_{j \in I \setminus \{k\}} f_j \right\} \right| \\
\notag &=  \sum_{i \in I} |\{f_i, f_k\}| + \sum_{i \in I} \left| \left\{ f_i, 1 - f_k \right\} \right| = 2  \sum_{i \in I} |\{f_i, f_k\}|  .
\end{align}
The first inequality in the proposition therefore follows from \Cref{conf_ess}. If there are $J \ge 1$ disjoint confined essential sets, say $U_1, \dots, U_J$, then $\cup_{i=1}^J U_i = \sqcup_{i=1}^J U_i$ and thus
\[ \int_M P_{\F} \, \omega \ge \int_{\cup_{i=1}^J U_i} P_{\F} \, \omega = \sum_{i=1}^J \int_{U_i} P_{\F} \, \omega \ge \sum_{i=1}^J 2 = 2J \, . \]

\QED
\vspace{12pt}

\begin{remark}
The proof of this result cannot be readily adapted to the more general case when $S_{\F} \ge 1$, although some intuition coming from the mean value theorem suggests that the result should remain valid.
\end{remark}\vspace{0pt}

%%%%%%%%%%%%%%%%%%%%%%%%%%%%%%%%%%%

\appendix

%%%%%%%%%%%%%%%%%%%%%%%%%%%%%%%%%%%%

%%%%%%%%%%%%%%%%%%%%%%%%%%%%%%%%%%%
\section{Displaceability in dimension two} \label{sec:disp_two}

An important fact which makes the Poisson bracket conjectures tractable in dimension two is the following characterization of displaceability for closed sets. Recall that in dimension two, a symplectic form is precisely an area form.\vspace{0pt}

\begin{proposition2}\label{charac_disp}
Let $(M, \omega)$ be a symplectic surface and $X \subset M$ be a closed connected subset. 
\begin{enumerate}
\item $X$ is displaceable if and only if $X$ is contained in a smoothly embedded closed disc $D \subset M$ of $\omega$-area less than half that of $M$.

\item If $X$ is displaceable,  its displacement energy is
\[ e_H(X) = \mathrm{inf} \, \left\{ \, \int_D \omega \, : \, X \subseteq D \subset M \mbox{ smoothly embedded closed disc }   \right\} \, . \]
\end{enumerate}
\end{proposition2} \vspace{6pt}

\noindent In its essence, this result appears to be well known to the experts. For instance, it is simply mentioned as a remark in \cite{BLT}, despite being of central importance in their arguments too.\footnote{In \cite{BLT}, a set is defined to be displaceable if its closure is displaceable in our sense. Therefore, as they remark, a connected set is displaceable in their sense if and only if it is contained in an open disc with smooth boundary of area strictly less than half the area of the surface. Our broader notion of displaceability, taken from \cite{PR}, considers open hemispheres in $S^2$ to be displaceable as well as some open discs whose boundaries are Osgood curves of large measure.} As we were not able to locate this characterization in the literature, we supply a proof of this characterization in this subsection. We first consider a particular case: \vspace{0pt}

\begin{lemma2}\label{disc_disj}
A smoothly embedded closed disc $D \subset (M, \omega)$ is displaceable if and only if its $\omega$-area is less than half the $\omega$-area of $M$, in which case its displacement energy equals its $\omega$-area.
\end{lemma2}\vspace{6pt}

\textit{Proof}. Set $c = \int_D \omega$ and $A = \int_M \omega$ (the latter being possibly infinite if $M$ is noncompact).

Assume $D$ is displaceable. Since $D$ is a compact subset of the Hausdorff space $M$, a small open neighbourhood of $D$ is also displaceable. By the area constraint, the area of this neighbourhood is at most $A/2$, hence $c < A/2$. That $e_H(D) \ge c$ follows from Usher's general and sharp energy-capacity inequality \cite{U}.

Conversely, assume $c < A/2$; we shall show that $D$ is displaceable and that $e_H(D) \le c$. Somewhat abusing notations, the smooth embedding $D \hookrightarrow M$ can be extended to a smooth embedding $D' \hookrightarrow M$ of a closed disc of area $c' \in (2c, A)$. As a consequence of Moser's trick, there are symplectic diffeomorphisms which are arbitrary close to mapping the pair $(D', D)$ to a pair of rectangles $(R', R)$ in $\R^2$ of corresponding areas, where $R$ lies inside a half of $R'$.  It is well known (see for instance the example in section 2.4 of \cite{P1}) that (small neighbourhoods of) $R$ can be displaced within $R'$ via a Hamiltonian isotopy of energy $c + \epsilon$ for any $\epsilon > 0$. By invariance of the energy under symplectomorphisms, we get $e_H(D) < c + \epsilon$ for every $\epsilon > 0$.

\QED\vspace{0pt}

\begin{remark2}\label{ham_isot}
We take the opportunity to reflect upon the proof of $e_H(D) \ge c$. Let $\phi_t : M \times [0,1] \to M$ be a Hamiltonian isotopy displacing $D$ and having Hofer energy $e$. We first consider the universal cover $M'$ of $M$ with the pullback symplectic form $\omega'$. It is well known (see \emph{e.g.} Corollary 1.8 in \cite{Ep}) that $M'$ is diffeomorphic to either $S^2$ or $\R^2$, and it is easy to construct global ‘action-angle’ coordinates turning this diffeomorphism into a symplectomorphism between $\omega'$ and the standard symplectic structures on $S^2$ or $\R^2$. Next, we note that the natural lift of $\phi_t$ to $M'$ is a Hamiltonian isotopy $\phi'_t$ with the same Hofer energy $e$ which displaces any given lift $D'$ of $D$. It follows that $e$ is at least the Hofer displacement energy of a disc of area $c$ inside $M'$, denote it $e_H(c;M')$. When $M' = \R^2$, Hofer's energy-capacity inequality \cite{H} states that $e_H(c;\R^2) \ge c$; this is also established in \cite{LM} via Gromov's nonsqueezing theorem. As this is true for any possible $e$, we have $e_H(D) \ge c$. If $M' = M = S^2$, things are more complicated. When the whole isotopy $\phi_t$ fixed a particular point $p \in S^2$, the displacement of $D$ effectively occurs in an open disc $S \subset \R^2$ of area $A$ and essentially the same argument as before applies. When no point is fixed by the isotopy, one can cook up a new Hamiltonian isotopy displacing $D$ and fixing a specific point at all times, but it is rather unclear to what extent the Hofer energy could increase in this way. Note that all of the above proofs of the energy-capacity inequality $e_H(D) \ge c$ use hard symplectic topology results, even though we are working only on surfaces; This can be traced to the fact that we need also consider the time-dependent Hamiltonian isotopies, making the problem implicitly higher dimensional. It would be much interesting to produce an ‘elementary’ proof of the energy-capacity inequality in the case of surfaces.
\end{remark2} \vspace{0pt}

In order to prove \Cref{charac_disp}, we shall need the following lemma. We give two proofs: the first one uses Lagrangian Floer theory and the second one uses more elementary and soft results. We thank Dominique Rathel-Fournier for discussions regarding the two approaches.\vspace{0pt}

\begin{lemma2} \label{circle_disj}
Let $(M, \omega)$ be a closed symplectic surface and $C \subset M$ a displaceable smoothly embedded circle. Then $C$ is contractible.
\end{lemma2}\vspace{0pt}

The case $M = S^2$ is trivial, so we shall assume $M \neq S^2$ below.\vspace{6pt}

\textit{Proof 1}. This proof uses Lagrangian Floer homology \cite{F}. Since $M \neq S^2$, $\pi_2(M) = 0$. We argue by contradiction: assume $C$ is noncontractible. Since this loop bounds no disc, $\pi_2(M,C) = 0$. Theorem $1$ in \cite{F} states that for all $\phi \in \Ham(M, \omega)$, $| C \cap \phi(C) | \ge cl_{\Z_2}(C) = 2$ where $cl_{\Z_2}(X)$ denotes the \emph{$\Z_2$-cuplength} of a topological space $X$, defined as the maximal integer $k$ such that there exist $k-1$ nonzero degree cohomology classes $\alpha_j \in H^*(X;\Z_2)$ satisfying $\alpha_1 \cup \dots \cup \alpha_{k-1} \neq 0$. Hence $C$ is not displaceable.

\QED
\vspace{6pt}

\textit{Proof 2}. Assume on the contrary that $C$ is noncontractible. Consider $\{\phi_t\}_{t \in [0,1]} \in \widetilde{\Ham}(M, \omega)$ displacing $C$ and set $C' = \phi_1(C)$. It is a classical fact due to Banyaga \cite{Ba} that Hamiltonian isotopies lie in the kernel of the \emph{flux morphism} (see also \cite{P1})
\[ \mathrm{Flux} : \widetilde{\mathrm{Symp}_0}(M, \omega) \to H^1_{dR}(M ; \R) : \{\psi_t\}_{t \in [0,1]} \mapsto  \,  \int_0^1 \left[ \frac{d \psi_t}{dt} \lrcorner \, \omega \right] \, dt \; , \]
from which we deduce that the isotopy  $\{\phi_t\}_{t \in [0,1]}$ generates a (usually degenerate) cylinder  $R : S^1 \times [0,1] \to M$ with boundary $-C + C'$ and area $\omega(R) = \int_{S^1 \times [0,1]} R^*\omega = 0$. Moreover, $C$ and $C'$ being embedded, isotopic, disjoint and noncontractible, it follows from Lemma 2.4 in \cite{Ep} that there exists an embedded cylinder $R' : S^1 \times [0,1] \hookrightarrow M$ such that $\partial R' = - C + C'$. Being embedded, its area $\omega(R') := \int_{S^1 \times [0,1]} (R')^*\omega$ satisfies $0 < \omega(R') < \omega(M)$. Since $R$ and $R'$ have the same boundary, these two $2$-chains differ by a $2$-cycle in $M$, hence the cohomology class $c = [R'] - [R] \in H_2(M ; \Z)$ has area $\omega(R') =  \omega(c) \in \omega(M) \Z$, which is a contradiction.

\QED
\vspace{6pt}

\begin{remark2}
It follows from \Cref{charac_disp} that \Cref{circle_disj} holds even without the assumption on $(M, \omega)$ being closed. Here we sketch a direct proof of this general version of the Lemma; It is a variation of the second proof above which has the additional merit of recasting the latter in more familiar terms. We thank the anonymous referee for suggesting this other argument.

Since Hamiltonian diffeomorphisms are compactly supported in the interior of $M$, $C$ necessarily belongs to the interior of $M$; We can therefore assume that $M$ is boundaryless. Suppose again that the circle $C$ is noncontractible in $M$, so that $\pi_2(M, C) = 0$; We shall prove that $C$ is nondisplaceable. From the homotopy long exact sequence associated to the pair $(M, C)$, we deduce that the map $ \pi_1(C) \to \pi_1(M)$ induced by the inclusion $C \subset M$ is injective. Let's denote its image $G$ and consider a covering map $p : M' \to M$ associated to $G$, so that $\pi_1(M') \cong G \cong \Z$. We shall identify $C$ with its preimage $p^{-1}(C)$. It follows from the classification of closed surfaces that $M'$ is not closed, hence open (since boundaryless). In fact, using the gradient flow of an exhausting Morse function on $M'$ without maxima (such a function exists since $M'$ is open, see \emph{e.g.} \cite{NRa}), $M'$ is seen to be diffeomorphic to a cylinder. Given an identification $M' \cong S^1 \times \R$, we can construct action-angle coordinates $(p, \theta)$ so that $C = \{p=0\}$ and $\omega' := p^*\omega$ reads $\omega' = \dd \lambda$ where $\lambda = p d \theta$. We thus reduced the problem to proving that the zero section $C$ of $(M' = T^*C, \omega' = \dd \lambda)$ is nondisplaceable, since any Hamiltonian diffeomorphism displacing $C$ in $M$ would lift to a Hamiltonian diffeomorphism displacing $C$ is $M'$. Assume on the contrary that there is a Hamiltonian diffeomorphism $\phi$ of $M'$ which displaces $C$. On the one hand, the flux $0 = \int_C \phi^*\lambda - \lambda = \int_{\phi(C)} \lambda$ is the $\omega'$-area of any chain between $C$ and $\phi(C)$. On the other hand, this $\omega'$-area is nonzero, since by \cite[Lemma 2.4]{Ep} there is such a chain which is an embedded cylinder.
\end{remark2}\vspace{6pt}

\textit{Proof of \Cref{charac_disp}}. Clearly, $X$ is displaceable whenever $X$ is contained in a displaceable set; In view of \Cref{disc_disj}, this is the case when $X$ is contained in a smoothly embedded closed disc $D$ with $\int_D \omega < (1/2) \, \int_M \omega$. In this case, the lemma also implies that
\[ e_H(X) \le \mathrm{inf} \, \left\{ \, \int_D \omega \, : \, X \subseteq D \subset M \mbox{ smoothly embedded closed disc }   \right\} \, .  \]

Now suppose that $X$ is displaceable. By definition, for every $\epsilon > 0$, there exists a compactly supported Hamiltonian isotopy $\phi_t$ with Hofer energy less than $e_H(X) + \epsilon$ displacing $X$ from itself. We shall prove that $X$ is contained in an embedded closed disc that is displaced by $\phi_1$ and which is thus included in the compact support of the isotopy $\phi_t$. This compact support is contained in the interior of a compact surface with finitely many boundary components, \emph{e.g.} a sublevel set associated to a sufficiently large regular value of a smooth exhaustion function defined on the interior of $M$. Since any compact surface with finitely many boundary components is symplectically embeddable in a closed symplectic surface (simply by capping the boundary circles by discs), we shall assume from now on that $(M, \omega)$ is closed. 

Since $X$ is assumed to be closed, we see that it is compact. Again, since $X$ is compact and $M$ is Hausdorff, the Hamiltonian isotopy displaces (the closure of) a small open neighbourhood $U$ of $X$. Let $\rho : M \to [0,1]$ be a smooth function such that $\rho^{-1}(0) = X$ and $\rho^{-1}(1) = M \setminus U$. By Sard's theorem, there exists a regular value $s \in (0,1)$ of $\rho$; The closed set $\rho^{-1}([0,s])$ has smooth boundary. Since $X$ is assumed connected, it is contained in some connected component $X_s$ of $X$. Of course, $X_s$ is displaced by $\phi_t$. 

\textsc{Claim}: \emph{For any sufficiently small $s \in (0,1)$, $X_s$ is contained in an embedded closed disc $D_s$ that is displaced by $\phi_t$. In fact, $X_s$ equals such a disc with finitely many holes in it. }

Assuming this claim, it follows from \Cref{disc_disj} that
 \[e_H(X) \le e_H(X_s) \le e_H(D_s) = \int_{D_s}  \omega < \;  \frac{1}{2} \, \int_M \omega \; ,\]
proving part (1) of the Proposition. We also have
\begin{align}
\notag e_H(X) + \epsilon \, &\ge \|\phi_t\|_H \ge e_H(D_s) = \int_{D_s} \omega \\
\notag & \ge  \mathrm{inf} \, \left\{ \, \int_D \omega \, : \, X \subseteq D \subset M \mbox{ smoothly embedded closed disc }   \right\} \, .
\end{align}
Taking the limit $\epsilon \to 0$ proves part (2) of the Proposition.

\QED \vspace{12pt}

\textit{Proof of the \textsc{Claim}}. We first recall some preliminary facts about curves in surfaces that can be proved via standard arguments along the lines of \cite[\S 1]{Ep}. The complement of a simple contractible smooth curve in an orientable connected surface consists of two connected components, one of which being a topological disc. These components are both topological discs if and only if the surface is the sphere. Given a smooth embedding of a closed disc inside the surface, its interior coincides with a connected components of the complement of its boundary circle.

We continue where we left in the proof of \Cref{charac_disp}. Since $M$ is compact, $\partial X_s \subseteq \rho^{-1}(s)$ and $s$ is a regular value of $\rho$, it follows that $\partial X_s$ consists of a finite number of disjoint smoothly embedded circles. From \Cref{circle_disj}, those boundary circles are contractible in $M$.  Moreover, since $X_s$ is connected, it lies completely within (the closure of) one of these two connected components.

Now, for each boundary circle, we select an embedded closed disc bounded by that circle as follows. (We will not distinguish between two embeddings that differ only by reparametrization.) When $M = \Sigma_g$ with $g \ge 1$, $M$ is aspherical and so there is a unique embedded disc bounded by the boundary circle. When $M=S^2$, there are two such discs, \emph{i.e.} the closures of the two connected components of the complement of the boundary circle; Decreasing $s$ if necessary (recall that there are finitely many boundary circles), we select the unique disc that has area no greater -- in fact, smaller -- than half the area of $M$.

We claim that $X_s$ belongs to one of these discs. Let's suppose the contrary. In view of the previous paragraph, this assumption implies that $X_s$ is contained in the complement of the interior of each of the selected discs. Then the union of $X_s$ with these discs defines a smooth surface without boundary $X'_s$ as well as an immersion from $X'_s$ to $M$ (induced by inclusion on each part of $X'_s$), which is thus a covering map. This map has degree one -- since only points from $X_s \subset X'_s$ are mapped to points in $X_s \subset M$ -- and is thus a diffeomorphism. In other words, $X_s$ and the interiors of the selected discs form a partition of $M$. Now, the fact that $\phi_1$ displaces the connected set $X_s$ means that $\phi_1$ sends $X_s$ in the interior of one of the selected disc, say the closed disc $D$. Since $\phi_1(\partial D)$ is a contractible simple curve, its complement in $M$ has two connected components, and $\phi_1(D)$ coincides with (the closure of) one of these two. Moreover, since $\phi_1(\partial D) \subset \mathrm{int} \, D$, it also bounds an open disc $D' \subset \mathrm{int} \, D$. Hence we have either that $\phi_1(D) = \overline{D'} \subset \mathrm{int} \, D$ or that $M = \phi_1(D) \sqcup D' = \phi_1(D) \cup D$. The first possibility is excluded at once by the area-preserving property of $\phi_1$. The second possibility is only possible if $M = S^2$, but is also excluded, as the area-preserving property of $\phi_1$ implies $\int_D \omega \ge (1/2) \, \int_M \omega$, which contradicts the way $D$ was chosen in the previous paragraph. 

Consequently, $X_s$ is contained in one of the selected discs, which we denote $D_s$. We note that $D_s$ contains all boundary circles of $X_s$ (its own boundary being one of them) and that those circles bound embedded discs inside $D_s$ (of area less than half that of $M$ when $M$ is the sphere). These discs therefore satisfy the properties that uniquely determined the discs we selected above, they coincide and it follows that our selected discs are all contained in $D_s$. This proves that $X_s$ equals $D_s$ minus finitely many discs. To prove that $\phi_1$ displaces $D_s$, we observe again that $\phi_1$ sends $X_s$ inside one of the connected components of its complement, which are the interiors of the selected discs other than $D_s$ as well as the complement of $D_s$. Essentially the same argument as in the previous paragraph excludes the first possibilities, hence $X_s$ is sent into the complement of $D_s$. The disc $\phi_1(D_s)$ coincides with the closure of one of the connected components of $M \setminus \phi_1(\partial D_s)$. Moreover, since $\partial D_s \subset X_s$, $\phi_1(\partial D_s)$  lies in the surface $M \setminus D_s$ and therefore bounds a disc in there that must coincide with one of the components of $M \setminus \phi_1(\partial D_s)$.  This implies either $\phi_1(D_s) \subset  M \setminus D_s$ or $D_s \subset \phi_1(D_s)$. The second possibility is excluded by the area-preserving property of $\phi_1$. We therefore conclude that $\phi_1(D_s) \cap D_s = \varnothing$.

\QED

%%%%%%%%%%%%%%%%%%%%%%%%%%%%%%%%%%%%%%%%%%

 \section{Short survey of Polterovich's conjecture} \label{sec:history}
 
We now describe in somewhat more details how the Poisson bracket invariants and the Poisson bracket conjectures were introduced and some of the progress made on these conjectures.\vspace{12pt}
 
We recall that $pb(\F)$, while nonnegative, may vanish:  Given a smooth function $h : M \to \R$ and a positive collection $\mathcal{G} = \{ g_i \}_{i=1}^N$ on $h(M)$ subordinate to some open cover $\mathcal{V} = \{V_i\}_{i=1}^N$ of $h(M)$, then $\F := \{ f_i := g_i \circ h \}_{i=1}^N$ is a positive collection on $M$ subordinate to the open cover $\U := \{ U_i := h^{-1}(V_i) \}_{i=1}^N$ such that $pb(\F) = 0$ and hence $pb(\U) = 0$. \vspace{12pt}

In comparison, when $\U$ consists of displaceable open sets, the invariant $pb(\U)$ cannot vanish if it is realized, \textit{i.e.} if there exists $\F \prec \U$ such that $pb(\F) = pb(\U)$. This follows from (the contrapositive of) the nondisplaceable fiber theorem \cite{EP}, which states that if a function $\vec{\F} : M \to \R^N : x \mapsto (f_1(x), \dots, f_N(x))$ has components which all pairwise Poisson commute, then some preimage of $\vec{\F}$ is nondisplaceable in $(M, \omega)$, hence any open set $U_i$ containing this fiber is also nondisplaceable. \vspace{12pt}

In \cite{EPZ}, Entov, Polterovich and Zapolsky generalized the nondisplaceable fiber theorem in a more quantitative way by considering partitions of unity subordinate to open covers; This result was reformulated in \cite{P2} in terms of the $pb$ invariant. We state here this last formulation in a way closer to the formulation of our previous results and which can be deduced from the material in \cite{PR}: Whenever $\U$ consists of displaceable open sets,
\[
 + \infty > pb(\U)e_H(\U) \ge \frac{1}{8N^2} \label{eq:Po} \tag{B.1}
\]
where $N$ is the cardinality of the cover $\U$ and where we defined $e_H(\U) := \max_{i \in \{1, \dots, N\}} \, e_H(U_i)$. In particular, $pb(\U)$ cannot vanish if $\U$ consists of (finitely many) displaceable sets. The proofs of the aforementioned results are sophisticated, relying on the functional analytic apparatus of (symplectic) quasi-states on the function space $C^{\infty}(M)$ (equipped with the Poisson bracket) constructed from spectral invariants obtained using the Hamiltonian Floer homology and the quantum cohomology of the symplectic manifold $(M, \omega)$. \vspace{12pt}

In \cite{P3}, Polterovich established that if the cover $\U$ is further assumed to be ``regular'' and ``fine'', morally meaning that each open set $U_i$ can be displaced within a sufficiently ``localized'' neighbourhood of it with the aid of a Hamiltonian diffeomorphism of energy smaller than some prescribed value $\mathcal{E}$ , then there exists a constant $C > 0$ depending on what is considered ``sufficient'' above, but not on the cardinality $N$ of the cover, such that $pb(\U)\mathcal{E} \ge C$. Notice that $e_H(U)$ is smaller than $\mathcal{E}$, possibly much smaller. Nevertheless, based on this result and on his intuition that ``irregular'' covers tend to have a higher $pb$ invariant, Polterovich asked whether the \Cref{conj} could be true. \vspace{12pt}

To our knowledge, \Cref{eq:Po} is still the best result valid without further assumption on $(M, \omega)$ or on $\U$; It can however be refined in some circumstances.  Motivated by the aforementioned ``local'' result of Polterovich, Seyfaddini \cite{Se} and Ishikawa \cite{I} studied more closely how the values of spectral invariants depend on localised data. Under some  monotonicity assumptions on $(M, \omega)$ and restricting $\U$ to consist either of images of symplectic embeddings of balls or of convex domains, respectively, they proved inequalities of the form $pb(\U)e_H(\U) \ge C/D^2$ where $C > 0$ is a universal constant and
\[ D = D(\U) = \mathrm{max}_{1 \le i \le N} \, \sharp \{ \, j \in \langle 1, N \rangle \, : \, \bar{U}_i \cap \bar{U}_j \neq \varnothing \,  \} \]
is what they call the \emph{degree of $\U$}. In fact, much like we do in this paper, Seyfaddini and Ishikawa establish somewhat stronger inequalities involving the capacities of the open sets rather than their (greater) displacement energy, thereby deducing that the displacement assumption in Polterovich conjecture is at best a sufficient condition for the nonvanishing of the $pb$ invariant.\vspace{12pt}

By a clever use of the lower semicontinuity of the $C^0$-norm of Poisson bracket on pairs of functions (see for instance \cite{PR}), Polterovich \cite{P4} and Buhovsky and Tanny \cite{BT} established a close variant of the strong conjecture for a large class of covers $\U$. Namely, given a Riemannian metric $g$ on $M$ compatible with $\omega$, there exist $\epsilon_0, C > 0$ depending only on $(M, \omega, g)$ such that for any $0 < \epsilon < \epsilon_0$, if $\U$ consists of open sets $U_i$ with diameter less than $\epsilon$, then for any $\F \prec \U$, $pb(\F)\epsilon^2 > C$. We note that $e_H(\U) \le c \, \epsilon^2$ for some constant $0 < c = c(M,g,\omega)$. Buhovsky--Tanny moreover proved in the same paper that the Poisson bracket conjecture is sharp, in the sense that it is possible to exhibit a family of such covers $\{\U_j\}_{j \in \N}$ and a sequence of positive numbers $\{\epsilon_j\}_{j \in \N}$ such that each $\U_j$ consists of sets of diameters at most $\epsilon_j$ and $\lim_{j \to \infty} \epsilon_j^2 = 0$, but $pb(\U_j)\epsilon_j^2 < C'$ for some $C' < + \infty$ independent of $j \in \N$. \vspace{12pt}

The situation on surfaces is more tractable than for general symplectic manifolds. On the one hand, there is the elementary and well-known characterization of displaceability in dimension $2$ which we proved in \Cref{sec:disp_two} for completeness: A closed set $X \subset M$ is displaceable if and only if it is contained in a closed smoothly embedded disc of area at most half that of $M$. Combined with the behavior of $pb$ with respect to refinements of open covers, the validity of a Poisson bracket conjecture on surfaces is essentially reduced to its validity on open covers by displaceable discs.  On the other hand, by studying the $L^1$-norm of the Poisson bracket functions $P_{\F}$, Buhovsky--Tanny \cite{BT} obtained several better lower bounds on $pb$ valid uniformly on all surfaces; Their results come into two sets of estimates, which we respectively dub degree estimates (which involve the \emph{degree} of a cover) and essential estimates (which involve the existence of so-called \emph{essential} discs to the cover). Explicitly, they proved that there exists a constant $C > 0$ such that for any closed symplectic surface $(M, \omega)$  and any open cover $\U$ of $M$ made of displaceable open \emph{discs},
\begin{align}
\notag & pb(\U) e_H(\U) \ge C \max \left\{ \chi(\mathcal{J}) \, , \, \bar{D}^{-2}    \right\} \, , \\
\notag & pb(\U) \mathrm{Area}(M, \omega) \ge C \max \left\{ |\mathcal{J}| \, , \, (\log \bar{D})^{-1}   \right\} \; .
\end{align}
Here, $\bar{D} = \bar{D}(\U) := \mathrm{max}_{1 \le i \le N} \, \sharp \{ \, j \in \langle 1, N \rangle \, : \, U_i \cap U_j \neq \varnothing \,  \}$ is what Buhovsky and Tanny call the \emph{degree of $\U$}, $\mathcal{J} = \mathcal{J}(\U) \subseteq \U$ is the subset of essential sets of $\U$, where $U_i \in \U$ is \emph{essential} if $\U \setminus \{U_i\}$ is not a cover of $M$, $|\mathcal{J}|$ is the cardinality of $\mathcal{J}$, and $\chi(\mathcal{J}) = 1$ if $\mathcal{J} \neq \varnothing$ and $0$ otherwise. These estimates follow from elementary, yet clever (and for the degree estimates, at times intricate) arguments with a strong geometric flavour. \vspace{12pt}

In an updated version of \cite{BT}, Buhovsky, Logunov and Tanny \cite{BLT} proved the Poisson bracket conjecture for every closed symplectic surfaces and for a universal constant $C$ \textit{i.e.} independent from $(M, \omega)$. They in fact accomplished more: Given \emph{two} partitions of unity $\F = \{f_1, \dots, f_N\}$ and $\GG = \{g_1, \dots, g_L\}$ on $M$, the authors considered the function
\[ P_{\F, \GG} : M \to [0, \infty) : x \mapsto \sum_{i=1}^N \sum_{j=1}^L |\{f_i, g_j \}(x)|  \, .\]
The quantity $\|P_{\F, \GG}\|$ then generalizes $pb(\F)$, since $P_{\F} = P_{\F,\F}$. This sort of invariant (an instance of which was already considered in \cite{P3}) could be interpreted as a measure of Poisson non-commutativity of (or of Poisson interaction between) the two partitions of unity, so that $pb(\F) \simeq pb(\F, \F)$ becomes a measure of Poisson self-interaction. Buhovsky--Logunov--Tanny proved that for partitions of unity $\F$ and $\mathcal{G}$ respectively subordinate to open covers $\U = \{U_1, \dots, U_N\}$ and $\V = \{V_1, \dots, V_L\}$ of $(M, \omega)$ consisting of displaceable open sets,
\begin{equation} 
\notag \int_M P_{\F,\GG} \, \omega \, \ge \, \frac{\mathrm{Area}(M, \omega)}{2 \, \mathrm{max}\{ e_H(\U), e_H(\V) \}} \; ,  \end{equation}
which readily implies
\[ \|P_{\F, \GG}\| \, \mathrm{max}\{ e_H(\U), e_H(\V) \} \, \ge \, 1/2  \; .  \]
Loosely speaking, they achieved this by noticing that it is possible to bound $\int_M P_{\F,\GG} \, \omega$ from below in terms of the numbers of intersection points of the level sets of the functions from $\F$ and $\mathcal{G}$, and that these numbers are themselves universally bounded from below. For comparison with the methods of the present paper, it is worth mentioning that their proof of the inequality in the case of only one open cover $\U$ also requires to establish estimates on pairs of open covers. \vspace{12pt}

More recently, Shi and Lu \cite{SL} adapted the arguments in \cite{BLT} to find a sufficient and necessary condition for an open cover by (not necessarily displaceable) discs in general position $\U$ on any closed symplectic surface to have nonvanishing $pb$ invariant: No two discs from $\U$ should suffice to cover $M$.\footnote{We pointed out on multiple occasions that this condition is necessary. In the case of surfaces of genus $g \ge 1$, the condition is automatically satisfied, since such surfaces have Lusternik-Schnirelmann category $3$; \Cref{cor_pb_thm} therefore implicitly establishes the sufficiency of this condition when $g \ge 1$. The sufficiency of the condition is thus most interesting in the case of the sphere, for which it nicely generalizes and simplifies our $3$-localization or confined star assumptions.} They moreover proved that, when this condition is satisfied, the weak Poisson bracket conjecture is valid, with
\begin{equation}\label{ShiLu}
\notag \int_M P_{\F} \, \omega \, \ge \, 2 \; .
\end{equation}
More generally, in the case of two open covers $\U$ and $\V$, their methods establish the following more general fact: The infimum of $P_{\F, \GG}$ over positive collections $\F \prec \U$ and $\GG \prec \V$ is positive if and only if $M \not \subset U \cup V$ for every $U \in \U$ and $V \in \V$, in which case $\int_{M} P_{\F, \GG} \, \omega \ge 2$. \vspace{12pt}

The results of the present paper were obtained around the same time as those of \cite{BLT}. We largely rely on the way Poisson brackets, displacement energies and areas behave under pullbacks along symplectic covering maps. This suggests that a convenient way of thinking about the Poisson bracket conjectures may be in terms of Poisson morphisms between symplectic manifolds (called \emph{symplectic submersions} in \cite{Pa}). This idea has been explored in the author's PhD thesis \cite{Pa} as a way to approach the Poisson bracket conjectures in higher dimensions and met with some success, \emph{e.g.} it yielded another proof of Polterovich--Buhovsky--Tanny's result on the $pb$ invariant of metrically small open covers. Explaining these results in further details shall be the object of a forthcoming paper.

\vspace{12pt}

%
% ----------  B I B L I O G R A P H I E  ----------
%

\end{document}